\documentclass{siamltex}
\usepackage{amsmath,amssymb}
\usepackage[table]{xcolor}%
\usepackage{graphicx}
\usepackage{amsfonts}
\usepackage{tikz}
\usepackage{placeins}
\usepackage{subfigure}
\usepackage[table]{xcolor}
\usepackage{booktabs}
\usepackage{epstopdf}
\usepackage{epsfig,epsf,psfrag}
\usepackage{multirow}
\usepackage{mathtools}
\usepackage{latexsym,amssymb, amsfonts, setspace, geometry}
\newcommand{\dsum}{\displaystyle \sum}
\newcommand{\dint}{\displaystyle \int}
\newcommand{\lb}{L_{\mathcal{B}}}
\newcommand{\twopartdef}[4]
{
	\left\{
		\begin{array}{ll}
			#1 & \mbox{if } #2 \\
			#3 & \mbox{if } #4
		\end{array}
	\right.
}
\newtheorem{remark}{Remark}
\newtheorem{example}{Example}[section] 
\definecolor{mycyan}{rgb}{0,0.95,0.95}

\title{Fast, adaptive, high order accurate discretization of the 
Lippmann-Schwinger equation in two dimensions}

\author{Sivaram Ambikasaran,
Carlos Borges,
Lise-Marie Imbert-Gerard
\and Leslie Greengard
}

\begin{document}
\maketitle

\begin{abstract}
	We present a fast direct solver for two dimensional scattering problems,
where an incident wave impinges on a penetrable medium with compact support. 
We represent the scattered field using a volume potential whose kernel is the 
outgoing Green's function for the exterior domain.
Inserting this representation into the governing partial differential equation,
we obtain an integral equation of the Lippmann-Schwinger type. 
The principal contribution here is the development of an automatically adaptive,
high-order accurate discretization based on a quad tree data structure
which provides rapid access to arbitrary elements of the discretized system 
matrix.
This permits the straightforward application of state-of-the-art algorithms for
constructing compressed versions of the solution operator. 
These solvers typically require $O(N^{3/2})$ work, where $N$ denotes the 
number of degrees of freedom. 
We demonstrate the performance of the method for a variety of problems in both 
the low and high frequency regimes.
\end{abstract}


%%%%%%%%%%%%%%%%%%%%%%%%%%%%%%%%%%%%%%%%%%%
\section{Introduction}%%%%%%%%%%%%%%
%%%%%%%%%%%%%%%%%%%%%%%%%%%%%%%%%%%%%%%%%%%

The problem of acoustic or electromagnetic scattering from penetrable media
arises in a variety of applications, from 
medical imaging to remote sensing, nondestructive testing, sonar, 
radar and geophysics. 
In the two-dimensional setting, the governing partial differential equation is 
the time-harmonic Helmholtz equation
\begin{equation}
\Delta u(x)+\kappa^2(1+q(x))u(x)=0,\quad x\in\mathbb{R}^2
\label{eqn_time_harmonic_Helm}
\end{equation}
where $u(x)$ is the total field and $\kappa$ is the background wavenumber,
and the perturbation (or {\em contrast function}) $q(x)$ has compact support,
say in a domain $\Omega$. Note that for $q(x) > -1$ the solution is 
propagating, while for $q(x) < -1$ it is evanescent.
Using the standard language of scattering theory, the total field $u(x)$ can be 
expressed as the sum of a known incident field $u^{\text{inc}}(x)$ and an
unknown scattered field $u^{\text{scat}}(x)$. In order to be well-posed, 
the latter must satisfy the Sommerfeld radiation condition~\eqref{eqn_Sommerfeld}.

\begin{equation}
\lim_{r \to \infty} r^{1/2}\left(\dfrac{\partial u}{\partial r}^{\text{scat}} - i\kappa u^{\text{scat}}\right) = 0 \,\,\,\, r= \vert x \vert
\label{eqn_Sommerfeld}
\end{equation}
The incoming field is assumed to satisfy the homogeneous equation
\begin{equation}
\Delta u^{\text{inc}}(x)+\kappa^2u^{\text{inc}}(x)=0
\label{eqn_time_harmonic_Helm_inc}
\end{equation}
in some neighborhood $D$ containing $\Omega$.
From eqs.~\eqref{eqn_time_harmonic_Helm} and~\eqref{eqn_time_harmonic_Helm_inc}, it follows
that the unknown scattered field satisfies 
\begin{equation}
\Delta u^{\text{scat}}(x) + \kappa^2 q(x) u^{\text{scat}} = -k^2 q(x) u^{\text{inc}}(x) \, .
\label{eqn_time_harmonic_Helm_scat}
\end{equation}

While PDE-based methods can be used to discretize (\ref{eqn_time_harmonic_Helm_scat})
directly (see \cite{ey11,Barnett2015,zd14} and the references therein), we choose to 
represent the scattered field $u^{\text{scat}}$ using a volume potential
\begin{equation}
u^{\text{scat}}(x) = V[\psi](x) = \int_{\Omega} G_{\kappa}(x,y) \psi(y) dy ,
\label{eqn_volume_potential}
\end{equation}
where $\psi(y)$ is an unknown density and $G_{\kappa}(x,y)$ is the 
outgoing free space Green's function, given by
\begin{equation}
G_{\kappa}(x,y) = \dfrac{i}4 H_0^{(1)}(\kappa\vert x - y \vert) \, .
\label{eqn_free_space_Green}
\end{equation}

It is well-known that the operator $V$ in eq.~\eqref{eqn_volume_potential} is 
bounded as a map from $L^2(\Omega)$ to $H^2(D)$ and compact
as an operator acting on $L^2(\Omega)$ 
\cite{Costabel}. Further, it is straightforward to see that,
substituting the representation~\eqref{eqn_volume_potential} into 
eq.~\eqref{eqn_time_harmonic_Helm_scat}, we obtain the Lippmann-Schwinger 
integral equation
\begin{equation}
\psi(x) + \kappa^2 q(x) V[\psi](x) = f(x) \,
\label{eqn_main_eqn}
\end{equation}
where $f(x) = -\kappa^2q(x)u^{\text{inc}}(x)$. 
This is an invertible (resonance-free) Fredholm equation of the second kind 
with a weakly singular kernel.

\begin{remark}
The Lippmann-Schwinger equation is often used to 
denote an alternative formulation \cite{Colton}:
\begin{equation}
u(x) - \kappa^2 V[uq](x) = u^{\text{inc}}(x) \, ,
\label{classicalLS}
\end{equation}
derived by convolving the original Helmholtz equation 
(\ref{eqn_time_harmonic_Helm}) with the governing Green's 
function $G_{\kappa}(x,y)$.
The equation (\ref{eqn_main_eqn}) has two advantages.
First, one is often interested in gradients of the scattered field.
This can be done from (\ref{eqn_volume_potential}),
with full precision by quadrature.
Using (\ref{classicalLS}), one would need to numerically differentiate
the computed solution, with the attendant loss of accuracy.
The formulation (\ref{eqn_main_eqn}) is also slightly easier
to work with, since the contrast function $q(x)$ appears as a (left) 
diagonal multiplier of the volume integral operator, whereas it 
appears inside the volume integral in the classical scheme.
\end{remark}

\vspace{.2in}

The Lippmann-Schwinger equation poses several numerical challenges.
First, it leads to a large, dense linear system of equations for $\psi(x)$.
Second, it may involve a complicated contrast function $q(x)$, requiring adaptive
mesh refinement for effective resolution. Third,
it may be ill-conditioned due to multiple scattering once the contrast $q(x)$ 
and $\kappa$ are large. The literature in this area is substantial, and we do not seek
to review it here. Some relevant prior work on volume integral-based methods, 
fast solvers, and numerical scattering theory includes
\cite{Amar1983,siva,amirhossein2014fast,bebendorf2005hierarchical,borm2003hierarchical,borm2003introduction,Bruno2003,Bruno2004,Bruno2005,chandrasekaran2006fast1,chandrasekaran2006fast,chen2002fast,huang2006,corona2014n,duan2009,Ethridge2001,Barnett2015,greengard2009fast,hackbusch2001introduction,ho2012fast,hoying1,hoying2,Hohage2006,kong2011adaptive,Lanzara2004,Biros,martinsson2009fast,martinsson2013,martinsson2005fast,
Sandfort2010,vainikko2000fast}. 

Our goal in this paper is to develop a fast, adaptive, high-order accurate, 
direct solver, which shares features with many of the algorithms in the 
papers cited above.
We concentrate, in particular, on developing a framework that permits
the rapid computation of entries in the 
governing system matrix, once the adaptive data structure has been specified.
With this capability in hand, it is straightforward
to make use of modern hierarchical direct solvers.

There has been relatively little attention paid to the adaptive discretization
issue, and fast direct solvers for 
volume integral equations are typically implemented on either uniform Cartesian or  
polar grids.
This is a perfectly sensible approach, especially when first developing fast solvers,
whether based on separation of variables and FFTs or 
hierarchical direct methods
\cite{aguilar2002,Bruno2004,Bruno2005,chen2002fast,duan2009,Barnett2015,vainikko2000fast}.

In the next section, we describe a high-order adaptive approach to discretization
of the unknown density $\psi$ in (\ref{eqn_main_eqn}), followed by a detailed
explanation of how one can rapidly compute the matrix entries of the fully
discretized system. We then discuss a fast solver for the Lippmann-Schwinger
equation using the HODLR method of \cite{siva, amirhossein2014fast}
and present numerical examples illustrating the performance of the scheme.

\section{Discretization} \label{sec-disc}

We assume that we are given a square domain $D$ which contains the support of
a known contrast function $q(x)$. By inspection of the 
Lippmann-Schwinger equation~\eqref{eqn_main_eqn}, $D$ clearly contains the support
of our right-hand side $f(x) = -\kappa^2q(x)u^{\text{inc}}(x)$, and therefore
the unknown $\psi(x)$ as well. We subdivide $D$ using an adaptive quad-tree,
ensuring that $q(x)$ and $f(x)$ are well-resolved, and that the number of 
points per wavelength
is greater than or equal to a user-specified parameter in each leaf node,
on which we impose a $p \times p$ grid (based on a uniform grid for $p=4$ 
and based on a tensor product Chebyshev grid for $p>4$).
The unknowns are taken to be the values of $\psi$
at the $p \times p$ grids on every leaf node, and we will seek to enforce
the integral equation at the same nodes, corresponding to a Nystr\"{o}m 
discretization of the integral equation. If we let $\vec{\psi}$ denote the 
vector of function values at all leaf node grid points and we let 
$\vec{f}$ denote the 
vector of right-hand side values at those same leaf node grid points, then
the discrete version of 
eq.~\eqref{eqn_main_eqn} takes the form
\begin{equation}
	A \vec{\psi} = \vec{f}
	\label{eqn_linear_system} \, .
\end{equation}

\subsection{The adaptive data structure}
\label{section_data_structure}

\begin{figure}[htbp]
\begin{center}
\begin{tikzpicture}[scale=2]
\draw (0,0) rectangle (1,1);
\draw (0,0.5) -- (1,0.5);
\draw (0.5,0) -- (0.5,1);
\node at (0.5,0.5) {$\mathcal{B}$};
\node at (0.25,0.25) {$\mathcal{B}_0$};
\node at (0.75,0.25) {$\mathcal{B}_1$};
\node at (0.75,0.75) {$\mathcal{B}_2$};
\node at (0.25,0.75) {$\mathcal{B}_3$};
\end{tikzpicture}
\caption{Box $\mathcal{B}$ and its four children}
\label{refine_box}
\end{center}
\end{figure}
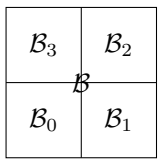

The domain $D$ is decomposed hierarchically, as follows.
Grid level 0 is defined to be the domain $D$ itself.
Grid level $l+1$ is obtained by recursive subdivision of each box 
$\mathcal{B}$ at level $l$ into four {\em child} boxes 
(Fig. \ref{refine_box}). $\mathcal{B}$ is referred
to as their {\em parent}. We allow for adaptivity, so that not all boxes at 
level $l$ are necessarily divided. We will, however, require that the tree 
satisfy a standard restriction - namely, that 
two leaf nodes which share a boundary point must be no more than one refinement 
level apart. Refinement is controlled by the following two criteria.
First, given the wavenumber (Helmholtz parameter) $\kappa$, a box is subdivided
if it is of side length $\lb$ and $\kappa \lb \geq 2 \pi M$, for a user-specified 
parameter $M$.  This ensures that there are at least $M$ points per
wavelength in each linear dimension.
Second, we ensure that the contrast $q(x)$ and the right hand side 
$f(x)$ are both \emph{well-resolved}.
For this, suppose we are given a leaf node $\mathcal{B}$ at level $l$.
We evaluate the 
functions $q(x)$ and $f(x)$ at tensor product Chebyshev nodes on $\mathcal{B}$
and compute the coefficients of the corresponding $p$th order Chebyshev 
approximation, which we will denote by $P_{\mathcal{B}}$.
We then subdivide the box into four child boxes $\mathcal{B}_i$, and compute 
$q(x)$ and $f(x)$ at tensor product Chebyshev nodes on each one.
If the computed values in the chid boxes agree with $P_{\mathcal{B}}$ to
a user-specified tolerance $\epsilon$, we terminate the refinement at box 
$\mathcal{B}$. Otherwise, we add the children to the data structure at level
$l+1$ and continue until the approximation criterion is satisfied.

\begin{figure}[htbp]
\begin{center}
\subfigure[Adaptive quad tree before level restriction]{
\begin{tikzpicture}[scale=1.5]
\draw (-2,-2) grid (2,2);
\draw (0.5,0) -- (0.5,1);
\draw (0,0.5) -- (1,0.5);
\draw (0,0.75) -- (0.5,0.75);
\draw (0.25,0.5) -- (0.25,1);
\end{tikzpicture}
}
\subfigure[Adaptive quad tree after level restriction]{
\begin{tikzpicture}[scale=1.5]
\draw (-2,-2) grid (2,2);
\draw (0.5,0) -- (0.5,1);
\draw (0,0.5) -- (1,0.5);
\draw (0,0.75) -- (0.5,0.75);
\draw (0.25,0.5) -- (0.25,1);
\draw (-0.5,0) -- (-0.5,2);
\draw (-1,0.5) -- (0,0.5);
\draw (-1,1.5) -- (1,1.5);
\draw (0.5,1) -- (0.5,2);
\end{tikzpicture}
}
\end{center}
\caption{An adaptive quad-tree and its level-restricted refinement.}
\label{fig_level_restricted}
\end{figure}
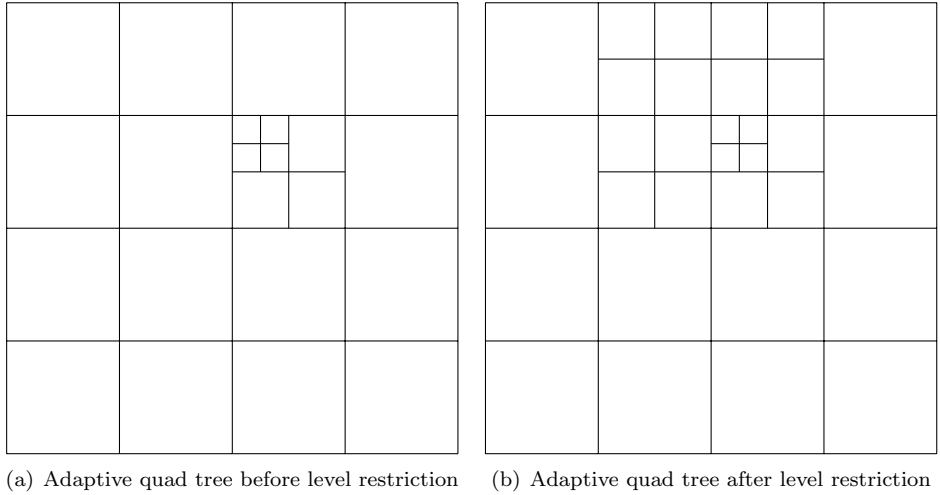

The procedure described above will not, in general, produce a
{\em level-restricted} tree. It is straightforward, however, to add further
refinements until that criterion is satisfied as well
(Fig.~\ref{fig_level_restricted}).
 
\subsection{The volume integral}
\label{section_volume_integral}

\begin{figure}[htbp]
\begin{center}
\includegraphics[scale=.4]{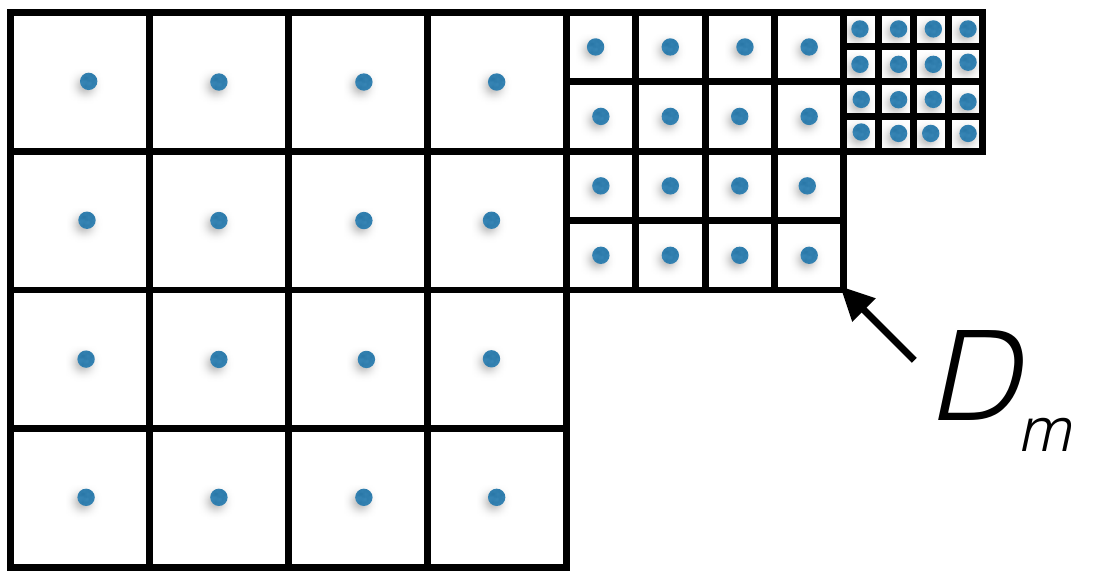}\label{figure_gridpts}
\end{center}
\caption{Discretization points for a leaf node $D_m$ and two neighboring
boxes at different levels of the hierarchy. A uniform grid is used for
4$th$ order accuracy and a tensor-product Chebyshev mesh for higher order
accuracy.}
\label{figure_points}
\end{figure}

We shall assume that there are $M$ leaf nodes in the tree, denoted by 
$\{D_m\}$, so that
\[ D = \cup_{m=1}^M D_m \, . \]
Thus, the volume integral defining the scattered field 
(\ref{eqn_volume_potential}) can be written in the form:
\begin{equation}
	u^{\text{scat}}(x) = V[\psi](x) = \dsum_{m=1}^M \dfrac{i}4 \int_{D_m} H_0^{(1)}(\kappa \vert x- y \vert) \psi(y)dy
	\label{eqn_volume_integral_leaf}
\end{equation}
Note that, since we are using an adaptive tree, 
the various $D_m$ can be at arbitrary levels in the spatial hierarchy. 
To obtain a Nystr\"{o}m method, it remains to discuss the computation 
of a high-order accurate discretization of $V[\psi]$. 
For this, we will require some notation. We let 
the values of the unknown density $\psi$ on leaf node $D_m$ be given 
by $\psi_{m,j} \approx \psi(x_j)$ for $j = 1,\dots,p^2$ where 
$x_j$ is one of the $p^2$ grid points on $D_m$ (Fig.~\ref{figure_gridpts}).

When clear from context, we write the full unknown vector as 
$\psi_j$ for $j=1,\dots,N$ where $N = p^2 M$ is the total number of unknowns.
We write $q_j = q(x_j)$ and $f_j = f(x_j)$ for the contrast function and
and the right-hand side, respectively,  at the corresponding grid points.

\subsubsection{Polynomial approximation}
\label{polynomial_approximation}

In order to obtain $p^{th}$ order accuracy, we build a $p$th order polynomial
approximation to the density $\psi$ on each leaf node $\mathcal{B}$  of length
$\lb$ centered at $(\mathcal{B}_1,\mathcal{B}_2 )$.
For this, we let $x = (\xi_1,\xi_2)$ and let $b_j(\xi_1,\xi_2)$ denote 
a suitable basis
for polynomials of two variables of degree $p-1$. 
That is, the $b_j(\xi_1,\xi_2)$ should
span $\{\xi_1^a\xi_2^b: 0 \leq a+b <p \text{ with }a,b \in \mathbb{N}^+\}$. 
The number of such basis functions is clearly $N_p = \dfrac{p(p+1)}2$.
In practice, we use the simple monomials
$\xi_1^a\xi_2^b$ for fourth order accuracy and 
Chebyshev polynomials of the form $T_a(\xi_1)T_b(\xi_2)$ for higher order
schemes, scaled to the unit box $[-1/2,1/2]^2$. 

Suppose now that we are given a $p$th order polynomial
defining $\psi(\xi_1,\xi_2)$ on $\mathcal{B}$:
\begin{equation}
\psi_{\mathcal{B}}(\xi_1,\xi_2) =
\dsum_{l=1}^{N_p} c_\mathcal{B}(l) b_l \left(\dfrac{\xi_1-\mathcal{B}_1}{\lb},
\dfrac{\xi_2-\mathcal{B}_2}{\lb} \right) \, .
	\label{fig_poly_psi}
\end{equation}
For the tensor product grid points $x_i = (\xi_{i,1},\xi_{i,2})$
lying in $\mathcal{B}$, we define the interpolation matrix 
$Q: \mathbb{R}^{N_p} \rightarrow \mathbb{R}^{p^2}$ with
entries $Q_{il}$ by
\[ Q_{il} = b_l \left(\dfrac{\xi_{i,1} -\mathcal{B}_1}{\lb},
\dfrac{\xi_{i,2}-\mathcal{B}_2}{\lb} \right) \, ,   \]
so that
\[ 
\psi_{\mathcal{B}}(\xi_{i,1},\xi_{i,2}) =
\dsum_{l=1}^{N_p} Q_{il} c_\mathcal{B}(l) \, .
\]

Note that in eq.~\eqref{fig_poly_psi}, the basis functions 
are independent of the level of box $\mathcal{B}$ in the quad tree, 
since we normalize by the box dimension $\lb$.
If we let $\vec{\psi_{\mathcal{B}}} \in \mathbb{R}^{p^2}$ denote the 
function values at the tensor product grid points in standard ordering, 
and $\vec{c}_{\mathcal{B}} =
(c_\mathcal{B}(1), c_\mathcal{B}(1),\dots, c_\mathcal{B}(N_p))
 \in \mathbb{R}^{N_p}$, then
$\vec{\psi_{\mathcal{B}}} \approx Q \vec{c_{\mathcal{B}}}$.
Let us denote by $Q^\dag$ the pseudoinverse of $Q$, so that the coefficients 
$\vec{c_{\mathcal{B}}}$ can be computed from $\vec{\psi_{\mathcal{B}}}$ via
\[ 
\vec{c_{\mathcal{B}}} = Q^\dag \vec{\psi_{\mathcal{B}}}  \, .
\]
The cost of computing $Q^\dag$ by means of the singular value
decomposition is negligible. Moreover, this is done only once and 
the pseudoinverse is stored.

Suppose now that we wish to compute an arbitrary matrix entry $V_{ij}$ 
corresponding to a $p$th order accurate Nystr\"{o}m discretization
of $V[\psi]$ in (\ref{eqn_volume_integral_leaf}).
We assume that $x_j$ lies in box $\mathcal{B}$ 
and that $x_i$ is an arbitrary point in the adaptive quad tree data structure
(including possible $\mathcal{B}$ itself).
Then, for $p$th order accuracy,
the contribution of the $j$th grid point to the $i$th target point is
precisely
\begin{equation}
V_{ij} \equiv \dfrac{i}4  \sum_{l=1}^{N_p} Q^\dag_{lj} \int_{\mathcal{B}} 
H_0^{(1)}(\kappa \vert x_i- y \vert) 
b_l \left(\dfrac{y_{1} -\mathcal{B}_1}{\lb},
\dfrac{y_{2}-\mathcal{B}_2}{\lb} \right) \, dy_1dy_2 \, , 
\label{Vij_def}
\end{equation}
where $y = (y_1,y_2)$.
This follows from the fact that
the $i$th column of 
$Q^\dag$ is a vector in $\mathbb{R}^{N_p}$, consisting of
the contributions (or projections) of the $i$th grid point 
onto each of the $N_p$ basis functions. 
Thus, the $(i,j)$ entry of the full matrix in the linear system in 
eq.~\eqref{eqn_linear_system} is given by
\begin{equation}
	A_{ij} = \delta_{ij} + \kappa^2 q(x_i) V_{ij} \, .
	\label{eqn_entry}
\end{equation}

It is worth emphasizing that the amount of work in 
computing $V_{ij}$ (or $A_{ij}$) is independent of 
$N$ (the total number of discretization points). Each entry could,
for example, be computed by adaptive quadrature on the fly.
We wish, however, to be able to evaluate an arbitrary matrix entry in 
only a few floating point operations. 
The remainder of this section is devoted to an extremely efficient
method for this task, based on the geometric relations 
between $x_i$ and $x_j$.

\begin{definition}
The
\textbf{colleagues} of a box $\mathcal{B}$ are boxes at the same level 
in the tree hierarchy which share a boundary point with $\mathcal{B}$. 
($\mathcal{B}$ is considered to be a colleague of itself.) 
Note that each box has at most $9$ colleagues (Fig.~\ref{fig_colleagues}). 
\begin{figure}[htbp]
	\begin{center}
	\begin{tikzpicture}
		\draw[white,fill=blue!30] (1,1) rectangle (2,2);
		\draw (0,0) grid (3,3);
		\node at (0.5,0.5) {$0$};
		\node at (1.5,0.5) {$1$};
		\node at (2.5,0.5) {$2$};
		\node at (0.5,1.5) {$3$};
		\node at (1.5,1.5) {$4$};
		\node at (2.5,1.5) {$5$};
		\node at (0.5,2.5) {$6$};
		\node at (1.5,2.5) {$7$};
		\node at (2.5,2.5) {$8$};
		\node at (1.5,1.775) {$\mathcal{B}$};
	\end{tikzpicture}
	\end{center}
\caption{Colleagues}
\label{fig_colleagues}
\end{figure}
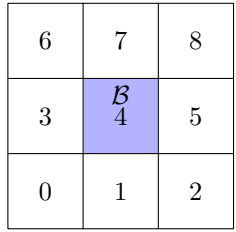

The \textbf{coarse neighbors} of $\mathcal{B}$ are leaf nodes that are one 
level higher than that of $\mathcal{B}$ and which share a boundary 
point with $\mathcal{B}$. Note that there can be at most $12$ coarse 
neighbors (Fig.~\ref{fig_coarse_neighbors}).
\begin{figure}[htbp]
	\begin{center}
	\subfigure[]{
	\begin{tikzpicture}[scale=0.9]
			\draw [fill=blue!30](0,0) rectangle (1,1);
			\draw (-2,-2) rectangle (0,0);
			\draw (1,-2) rectangle (3,0);
			\draw (1,1) rectangle (3,3);
			\draw (-2,1) rectangle (0,3);
			\node at (-1,-1) {$0$};
			\node at (2,-1) {$3$};
			\node at (2,2) {$6$};
			\node at (-1,2) {$9$};
			\node at (0.5,0.5) {$\mathcal{B}$};
	\end{tikzpicture}
	}
	\subfigure[]{
	\begin{tikzpicture}[scale=0.9]
			\draw [fill=blue!30](0,0) rectangle (1,1);
			\draw (-1,-2) rectangle (1,0);
			\draw (1,-1) rectangle (3,1);
			\draw (0,1) rectangle (2,3);
			\draw (-2,0) rectangle (0,2);
			\node at (0,-1) {$1$};
			\node at (2,0) {$4$};
			\node at (1,2) {$7$};
			\node at (-1,1) {$10$};
			\node at (0.5,0.5) {$\mathcal{B}$};
	\end{tikzpicture}
	}
	\subfigure[]{
	\begin{tikzpicture}[scale=0.9]
			\draw [fill=blue!30](0,0) rectangle (1,1);
			\draw (0,-2) rectangle (2,0);
			\draw (1,0) rectangle (3,2);
			\draw (-1,1) rectangle (1,3);
			\draw (-2,-1) rectangle (0,1);
			\node at (1,-1) {$2$};
			\node at (2,1) {$5$};
			\node at (0,2) {$8$};
			\node at (-1,0) {$11$};
			\node at (0.5,0.5) {$\mathcal{B}$};
	\end{tikzpicture}
	}
	\end{center}
\caption{Coarse neighbors}
\label{fig_coarse_neighbors}
\end{figure}
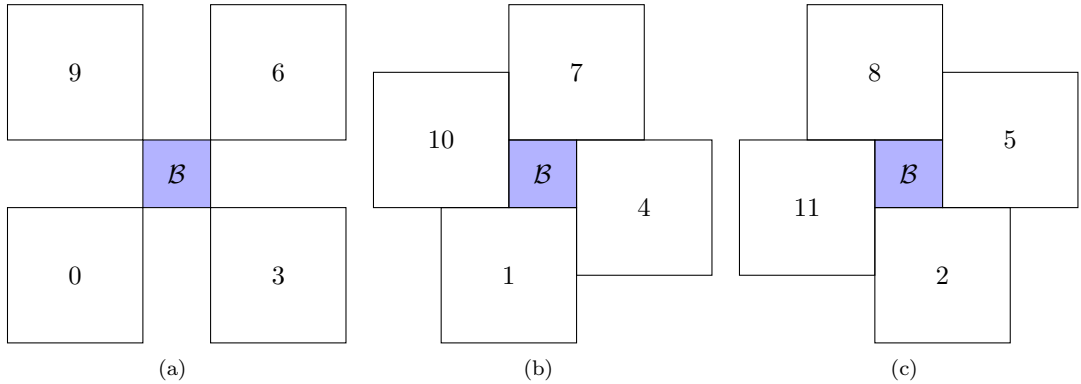

The \textbf{fine neighbors} of $\mathcal{B}$ are leaf nodes that are one 
level below that of $\mathcal{B}$ and share a boundary point with $\mathcal{B}$.
Note that there can be at most $12$ fine neighbors 
(Fig.~\ref{fig_fine_neighbors}).
\begin{figure}[htbp]
	\begin{center}
	\begin{tikzpicture}
		\draw (0,0) grid (4,4);
		\draw [fill=blue!30] (1,1) rectangle (3,3);
		\node at (0.5,0.5) {$0$};
		\node at (1.5,0.5) {$1$};
		\node at (2.5,0.5) {$2$};
		\node at (3.5,0.5) {$3$};
		\node at (3.5,1.5) {$4$};
		\node at (3.5,2.5) {$5$};
		\node at (3.5,3.5) {$6$};
		\node at (2.5,3.5) {$7$};
		\node at (1.5,3.5) {$8$};
		\node at (0.5,3.5) {$9$};
		\node at (0.5,2.5) {$10$};
		\node at (0.5,1.5) {$11$};
		\node at (2,2) {$\mathcal{B}$};
	\end{tikzpicture}
	\end{center}
\caption{Fine neighbors}
\label{fig_fine_neighbors}
\end{figure}

The \textbf{separated fine neighbors} of $\mathcal{B}$ are non-neighboring 
leaf nodes that are one level below that of $\mathcal{B}$ and 
whose parent shares a boundary point with $\mathcal{B}$. 
There are at most $16$ separated fine neighbors 
(Fig.~\ref{fig_intermediates}).
\begin{figure}[htbp]
	\begin{center}
	\begin{tikzpicture}
		\draw (0,0) grid (6,6);
			\draw [fill=white] (1,1) rectangle (5,5);
			\draw [fill=blue!30] (2,2) rectangle (4,4);
			\foreach \i in {0,...,4} {
				\node at (\i+0.5,0.5) {\i};
			}
			\foreach \i in {5,...,9} {
				\node at (5.5,\i-4.5) {\i};
			}
			\foreach \i in {10,...,14} {
				\node at (15.5-\i,5.5) {\i};
			}
			\foreach \i in {15,...,19} {
				\node at (0.5,20.5-\i) {\i};
			}
			\node at (3,3) {$\mathcal{B}$};
	\end{tikzpicture}
	\end{center}
\caption{Separated fine neighbors}
\label{fig_intermediates}
\end{figure}

All other leaf nodes are separated from $\mathcal{B}$ by at least the length
$\lb$ and are considered its 
\textbf{far field}.
\end{definition}

\subsection{Far-field interactions} 

The easiest case to consider is when $x_i$ lies in the {\em far field}
of box $\mathcal{B}$. If we let 
$x_i = (\xi_{i,1},\xi_{i,2})$ and denote by $\theta_i$ the polar
angle of $x_i$ with respect to the box center 
$(\mathcal{B}_1,\mathcal{B}_2)$, the 
Graf addition theorem \cite{abramowitz,handbook} states that
\begin{equation}
	H_0^{(1)}(\kappa |x_i-y|) = \dsum_{n=-\infty}^{\infty} 
H_n(\kappa \Vert (\xi_{i,1} - \mathcal{B}_1,\xi_{i,2} - \mathcal{B}_2) \Vert) 
J_n(\kappa \Vert (y_1 - \mathcal{B}_1,y_2-\mathcal{B}_2) \Vert) 
e^{in(\theta_i-\theta_y)}
	\label{eqn_Graf}
\end{equation}
where $\theta_y$ denotes the polar angle of the point $y \in \mathcal{B}$ 
with respect to the box center.
Since the dimensions of the leaf node are assumed to satisfy 
$\kappa \lb \leq 8$, we can truncate eq.~\eqref{eqn_Graf} after 
about $L+8$ terms with an error of approximately $2^{-L}$. It is easy to 
verify that
\begin{equation}
V_{ij} \approx \frac{i}{4} \sum_{n=-L}^{L} 
M_n(j) H_n(\kappa \Vert |x_i - (\mathcal{B}_1,\mathcal{B}_2)\Vert)
e^{il \theta_i}
\label{mpeval}
\end{equation}
where
\begin{equation}
M_n(j) = \sum_{l=1}^{N_p} C_{nl} Q^\dag_{lj}
\end{equation}
and
\begin{equation}
C_{nl} = 
\dint_{\mathcal{B}} 
J_n(\kappa \Vert (y_1 - \mathcal{B}_1,y_2-\mathcal{B}_2) \Vert) 
e^{-il\theta_y} 
b_l \left(\dfrac{y_1-\mathcal{B}_1}{\lb},\dfrac{y_2-\mathcal{B}_2}{\lb} \right) dy_1 dy_2 \, .
	\label{eqn_main_integral_scaled_truncated_2}
\end{equation}
Note that the integrals in 
(\ref{eqn_main_integral_scaled_truncated_2}) are smooth and need to be 
computed only once per level at negligible cost.
Moreover, $M_n(j)$ also only needs to be computed
once per level. Thus, the cost of computing a far field matrix entry
to fourteen digits of accuracy is essentially that of evaluating the 
multipole expansion (\ref{mpeval}) with $L = 45$, 
which can be done about one million times per second
on a single core.

\subsection{Separated fine neighbors}

The multipole expansion for $\mathcal{B}$ in the previous section is 
slowly converging for the separated fine neighbors (Fig. 
\ref{fig_intermediates}). However, it is straightforward to create four
child boxes from $\mathcal{B}$ and to compute
the multipole expansions for each of the four children from each basis
function 
$b_l \left(\dfrac{y_1-\mathcal{B}_1}{\lb},\dfrac{y_2-\mathcal{B}_2}{\lb} 
\right)$. Evaluating the corresponding multipole expansions 
at each separated fine
neighbor target point yields the desired value $V_{ij}$. 

\subsection{Near field interactions}

For points in the near field, we rewrite 
(\ref{Vij_def}) in the form:
\begin{equation}
V_{ij} \equiv  \sum_{l=1}^{N_p} Q^\dag_{lj}  V^\kappa_{il}  \, ,
\end{equation}
where
\begin{equation}
V^{\kappa}_{il} = \dfrac{i}4 
\dint_{\mathcal{B}} H_0^{(1)}(\kappa \vert x_i-(y_1,y_2) \vert) 
b_l \left(\dfrac{y_1-\mathcal{B}_1}{\lb},
\dfrac{y_{2}-\mathcal{B}_2}{\lb} \right) dy_1 dy_2 \, .
\label{eqn_main_integral}
\end{equation}

It will be convenient to rescale variables so that integrals are always
carried out over the unit box. Thus, instead of 
(\ref{eqn_main_integral}), we write
\begin{equation}
	V^{\kappa}_{il} = \dfrac{i\lb^2}4 
\dint_{U} 
H_0^{(1)}(\kappa \lb r_i) b_l \left(y_1,y_2 \right) dy_1 dy_2
	\label{eqn_main_integral_scaled}
\end{equation}
where $U = [-0.5,0.5]^2$ and $r_i=\vert x_i/\lb-(y_1,y_2) \vert$. 
Note that, in this representation, the integral 
$V^{\kappa}_{il}$ depends only on the relative coordinates of the target
point $x_i$ to the unit box $U$. There are only
finitely many such possibilities, corresponding to the various target locations
in the colleagues, coarse neighbors, or fine neighbors. 
Thus, the number of such interactions is fixed.
$x_i$ must be one of $p^2$ possible locations in each of $9$ possible 
colleagues, $12$ possible fine neighbors and $12$ possible coarse neighbors. 
Thus, we seek to compute tables of dimension $p^2 \times 9 \times p^2$,
$p^2 \times 12 \times p^2$, and 
$p^2 \times 12 \times p^2$, respectively.
For this, we use the fact \cite{abramowitz,handbook} that 
\[ H_0^{(1)}(z) = J_0(z) + iY_0(z) \]
with  
\[ J_0(z) = \sum_{m=0}^\infty \frac{(-1)^m}{(m!)^2} 
\left( \frac{z}{2} \right)^{2m}, 
\]
\[
Y_0(z) = \frac{2}{\pi} \left( \ln \left(\frac{z}{2} \right) + \gamma \right)
J_0(z) + \frac{2}{\pi} \left( 
\frac{\frac{1}{4}z^2}{(1!)^2} - 
(1 + \frac{1}{2}) \frac{(\frac{1}{4}z^2)^2}{(2!)^2} + 
(1 + \frac{1}{2} + \frac{1}{3}) \frac{(\frac{1}{4}z^2)^3}{(3!)^2} - \dots 
\right),
\]
where $\gamma$ is Euler's constant.
This permits us to write
\begin{equation}
H_0^{(1)}(\kappa \lb r) = 
\dsum_{p=0}^{\infty} c_p(\kappa \lb) \left(\dfrac{r}2\right)^{2p} + 
\dsum_{p=0}^{\infty}d_p(\kappa \lb) \left(\dfrac{r}2\right)^{2p} 
\log\left(\dfrac{r}2\right)
\label{eqn_Hankel_Bessel}
\end{equation}
where 
\[ c_p(z) = a_p(z) + \dfrac{2i}{\pi}g_p(z),\ 
d_p (z) = \dfrac{2i}{\pi} a_p(z),\ 
a_p(z) = \dfrac{(-z^2)^p}{(p!)^2}, \ 
g_p(z) = \left( \gamma + \log(z) - H_p \right)a_p(\kappa \lb). 
\]
Here, $H_0 = 0, H_1 = 1, H_2 =  (1 + \dfrac{1}{2}), 
H_3 = (1 + \dfrac{1}{2} + \dfrac{1}{3})$, etc. 
This permits us to write
\begin{align}
	V^{\kappa}_{il} & = 
\dfrac{i\lb^2}4 \dsum_{p=0}^{p_{\max}} c_p(\kappa \lb) \dint_{U} 
\left(\dfrac{r}2\right)^{2p} b_l \left(y_1,y_2 \right) dy_1 dy_2 \nonumber\\
& + \dfrac{i\lb^2}4 \dsum_{p=0}^{p_{\max}} d_p(\kappa \lb) 
\dint_{U} \left(\dfrac{r}2\right)^{2p} \log\left(\dfrac{r}2\right) 
b_l \left(y_1,y_2 \right) dy_1 dy_2 \, ,
\label{eqn_main_integral_scaled_truncated}
\end{align}
for a suitable choice of $p_{max}$.

The integrals over $U$ in eq.~\eqref{eqn_main_integral_scaled_truncated} are
independent of the wavenumber $\kappa$, so we may tabulate
these values to double precision accuracy using adaptive Gaussian quadrature. 
For $\kappa \lb \leq 8$, it is easy to verify that 
$p_{max} = 60$ is sufficient to obtain fourteen digits of accuracy.
The quantities $V^{\kappa}_{il}$ are then obtained using $O(p_{max})$
floating point operations. The storage requirements for this scheme are:
\begin{center}
\begin{tabular}{ll}
	$p^2 \times N_p \times 9 \times (2p_{\max}+2)$ & for colleagues\\
	$p^2 \times N_p \times 12 \times (2p_{\max}+2)$ & for fine neighbors\\
	$p^2 \times N_p \times 12 \times (2p_{\max}+2)$ & for coarse neighbors.
\end{tabular}
\end{center}

\vspace{.1in}

Unfortunately, for $\kappa \lb > 4$, the computation 
in \eqref{eqn_main_integral_scaled_truncated} can be subject to catastrophic
cancellation. (The series oscillates in sign and involves intermediate terms 
of large magnitude.) This loss of precision is easily overcome by 
subdividing $\mathcal{B}$ into four children and generating tables for each
child box. 
With this scheme, the storage costs remain negligible and the matrix entries 
$V_{ij}$ corresponding to local interactions are computed using only a few 
hundred floating point operations.

\section{Fast direct solvers}
\label{section_fast_direct}

The reason we have focused in this paper 
on the rapid computation of matrix entries
is that a new generation of fast solvers has been developed over the last 
decade or so which permits the inversion of equations such as the 
Lippmann-Schwinger equation in $O(N^{3/2})$ work. All that is required is 
a suitable ordering of the unknowns and access to matrix entries on the fly.
For the sake of simplicity, we have chosen to use the 
HODLR scheme of \cite{siva, amirhossein2014fast}.
This solver relies on 
a hierarchical partitioning of the matrix into a sequence of off-diagonal 
low-rank blocks (from which the acronym is derived). 
More precisely, the off-diagonal blocks are approximated using 
low-rank factorizations to a user-specified precision $\epsilon$.
We refer the reader to the original papers for details and to 
\cite{bebendorf2005hierarchical,borm2003hierarchical,borm2003introduction,
chandrasekaran2006fast1,chandrasekaran2006fast,chen2002fast,corona2014n,
ho2012fast,hoying1,goreinov1997theory,greengard2009fast,
hackbusch2001introduction,kong2011adaptive,liberty2007randomized,
martinsson2005fast,rjasanow2002adaptive,zhao2005adaptive} for 
related schemes. 
All of these fast solvers, of course, require sampling at most $O(N^{3/2})$
matrix entries.
For PDE-based analogues, see
\cite{Barnett2015,martinsson2009fast,hoying2} and the references
therein.

Several of the methods above have smaller constants than HODLR for volume
integral equations with singular kernels, but require a slightly more
complicated interface. In fact,
in the low frequency regime, some of the methods above require only 
$O(N \log N)$ work (see, for example, \cite{corona2014n,hoying2}). 
We postpone such accelerations to a later date.

\section{Numerical results}
\label{section_results}
In this section, we illustrate the performance
of our adaptive solver for the Lippmann-Schwinger equation.
In each example, we use an incoming plane 
wave propagating in the $x$-direction of the form
$u^{\text{inc}}=e^{ikx}$ and calculate
the value of the scattered field $u^{\text{scat}}$
on a square $D$ which contains the support of the contrast 
function $q(x)$.

We will make use of both uniform and adaptive grids.
As indicated in section \ref{sec-disc}, in the adaptive case,
we refine the domain $D$ using a quad-tree, halting refinement
of a leaf node when it is determined to have resolved the contrast
function $q(x)$ and the right-hand side $f(x)$ to a user-specified
tolerance $\epsilon$. In order to be conservative, 
we then use a tolerance of $\epsilon \times 10^{-4}$ in the HODLR solver.
The number of terms used in various far field expansions is
chosen as in \cite{crutchfield2006remarks}.
All simulations in this paper were carried out using a single core
of a 2.5GHz Xeon processor. 

\begin{example}Radially symmetric contrast functions \end{example}

In our first two examples, we consider a radially symmetric
contrast function $q(x)$, from which it is straightforward to compute the
exact solution using separation of variables, the fast Fourier transform,
and a high-order accurate ODE solver (see, for example, \cite{Barnett2015}). 
Setting the wavenumber $\kappa = 40$, we will use either a Gaussian
with $q(x) = 1.5e^{-160|x|^2}$ or a ``flat bump" with 
$q(x) = 0.5 \text{erfc}(5(|x|^2-1))$,
where $\text{erfc}$ is the complementary error function 
\cite{abramowitz,handbook}
\begin{equation*}
\text{erfc}(x)=\frac{2}{\sqrt{\pi}}\int_x^\infty e^{-t^2}~dt.
\end{equation*}

\begin{figure}[htbp]
\centering
\subfigure[Contrast function]{
\includegraphics[scale=.25]{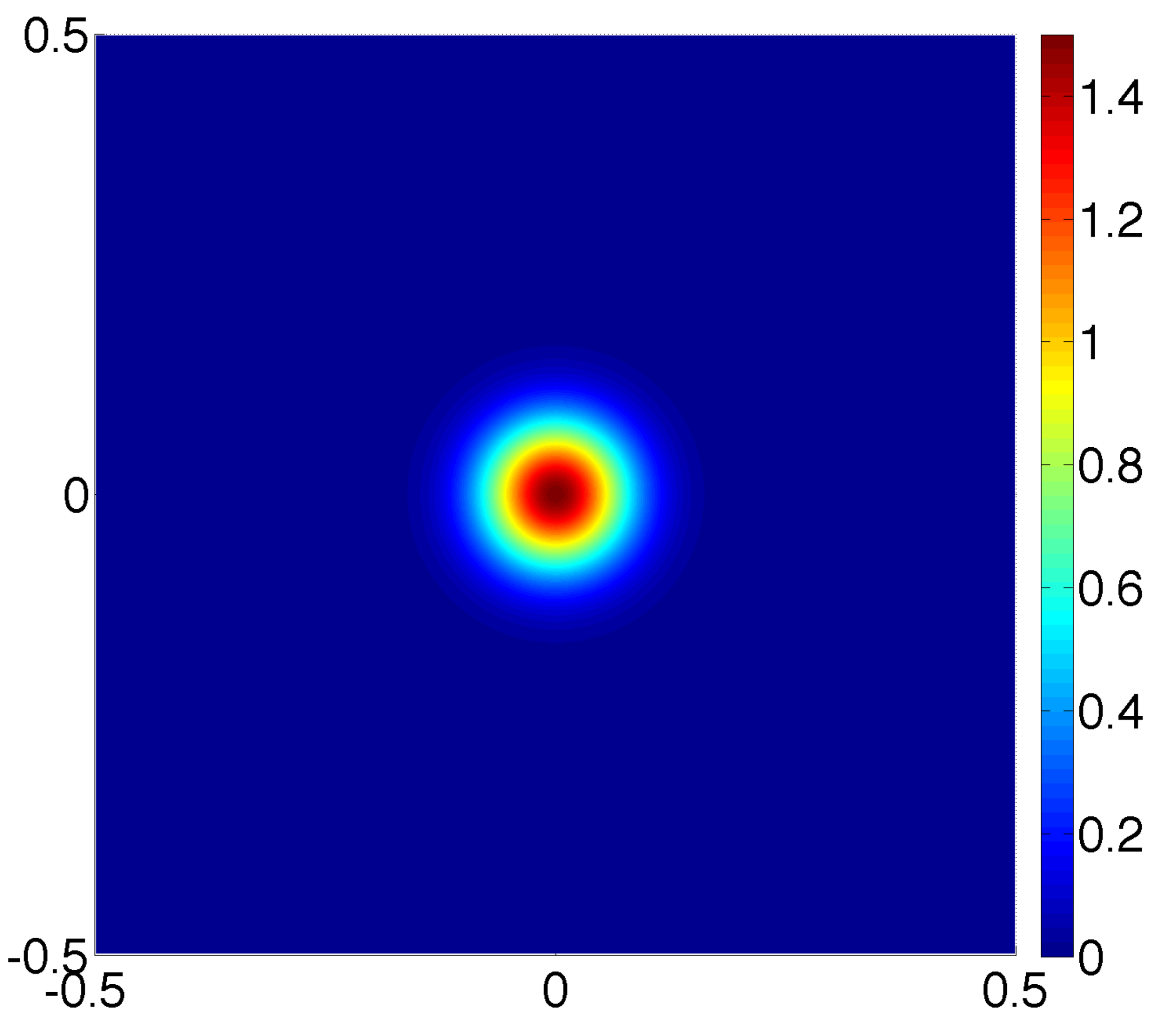}\label{figure_1_0}
}
\subfigure[Real part of total field]{
\includegraphics[scale=.25]{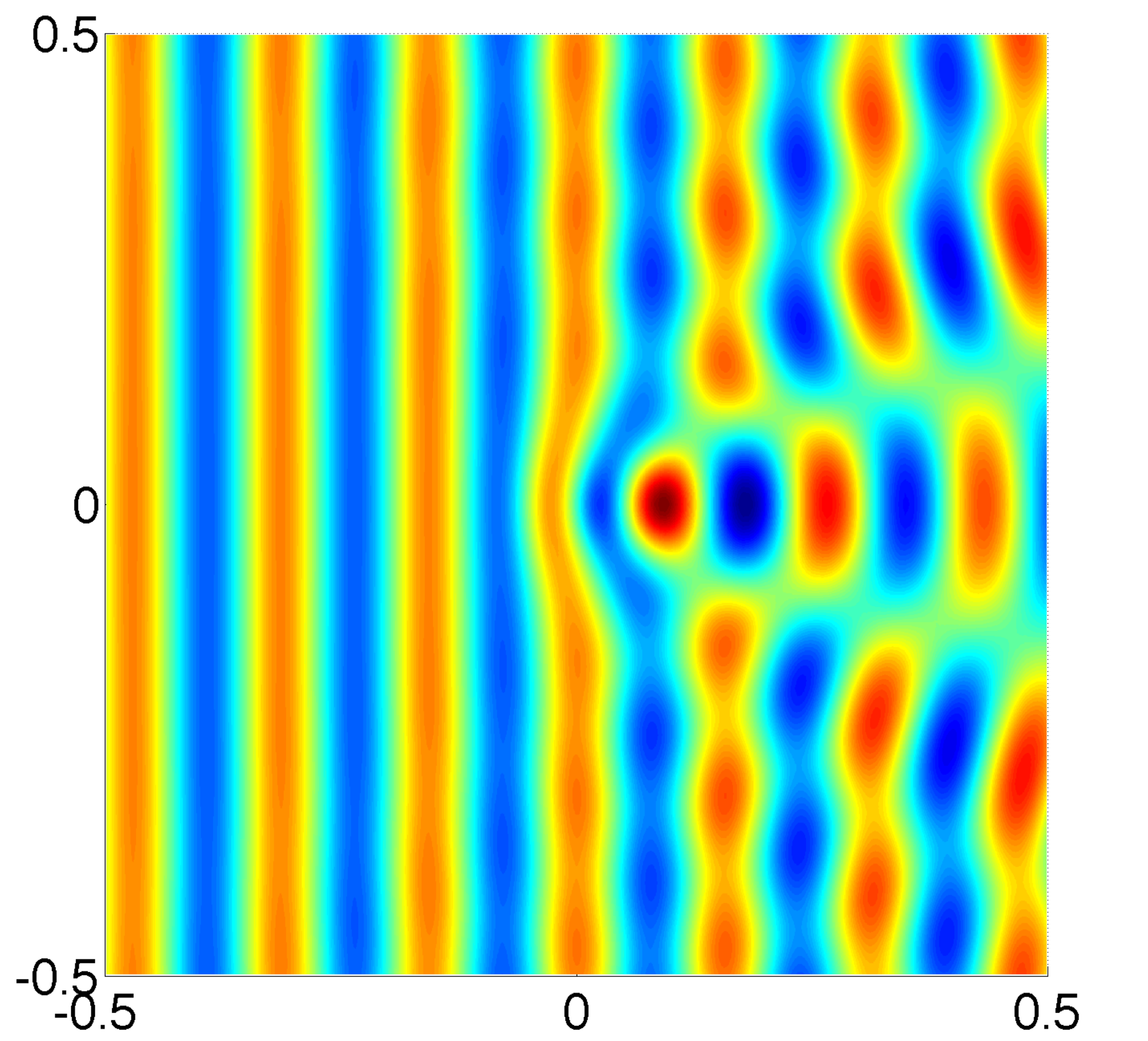}\label{figure_1_1}
}
%\caption{A discretization for the contrast function $q$ in the domain $D$.}
\label{figure_example_1}
\caption{Gaussian contrast}
\end{figure}

\begin{figure}[htbp]
\centering
\subfigure[Contrast function]{
\includegraphics[scale=.25]{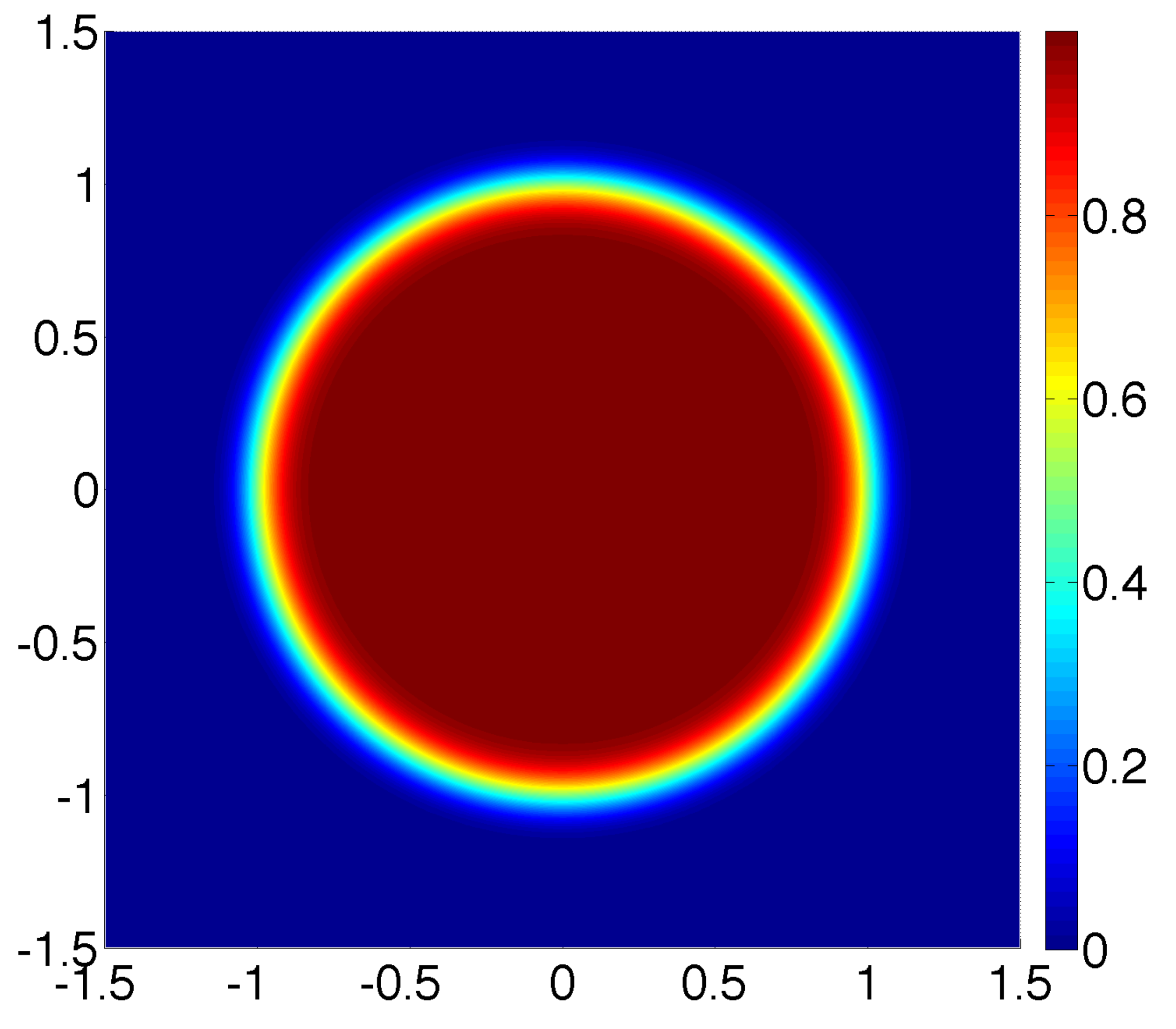}\label{figure_11_0}
}
\subfigure[Real part of total field]{
\includegraphics[scale=.25]{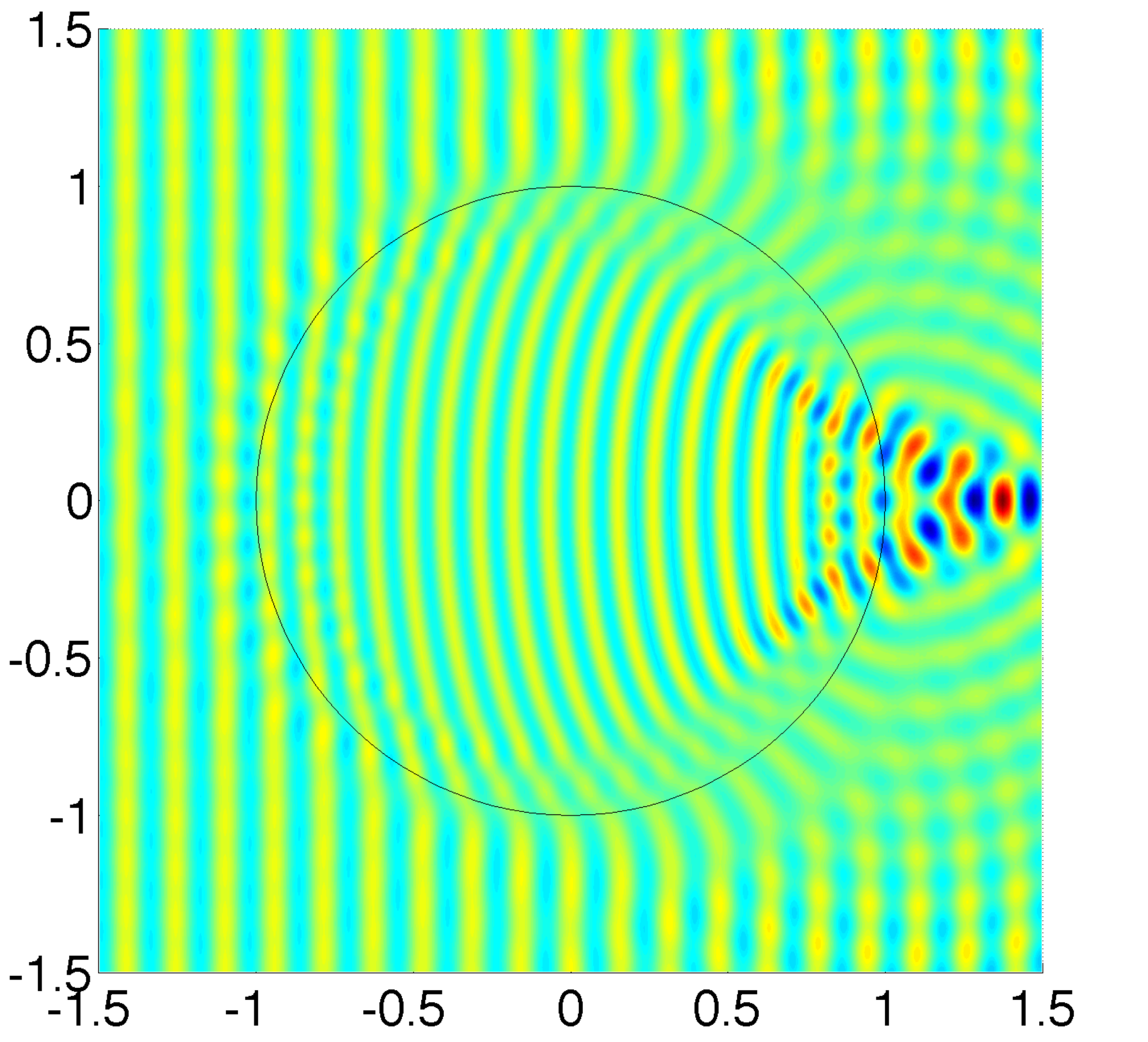}\label{figure_11_1}
}
\subfigure[Grid for adaptive solver]{
\includegraphics[scale=1.75]{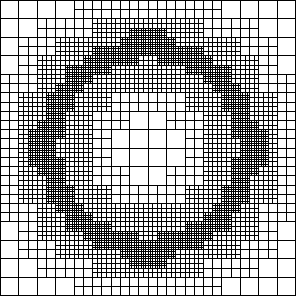}\label{figure_11_2}
}
\subfigure[Surface plot of contrast function]{
\includegraphics[scale=1.0]{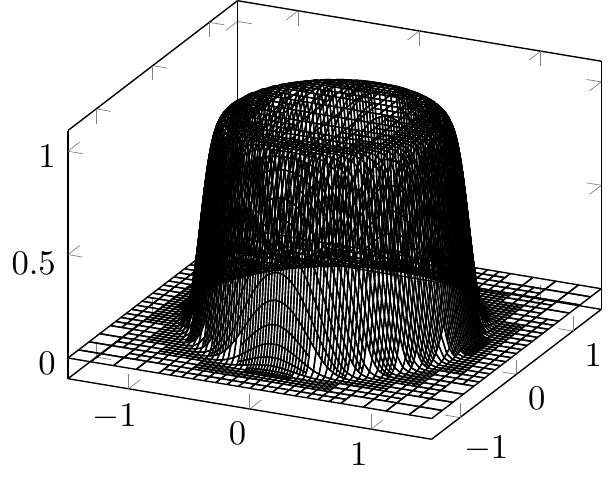}\label{figure_11_3}
}
\label{figure_example_11}
\caption{``Flat bump" contrast}
\end{figure}

In Fig. \ref{figure_1_0} we plot the Gaussian contrast, and in 
Fig. \ref{figure_1_1} we plot the real part of the total field 
after solving the scattering problem.
In Fig. \ref{figure_11_0} we plot contours of 
the ``flat bump" contrast function,
in Fig. \ref{figure_11_1} we plot the real part of the total field, 
in Fig. \ref{figure_11_2} we show the the adaptive discretization of 
the domain $D$, and in Fig. \ref{figure_11_3} we plot the contrast function 
$q(x)$. 

We let $\tilde{u}$ denote the total field calculated using our solver, 
we let $u$ denote the total field computed using separation of variables
as in \cite{Barnett2015}. 
Table \ref{gauss_uniform} shows the performance
of the solver on a uniform grid.

\begin{table}[h]
\begin{center}
\caption{Convergence behavior of Lippmann-Schwinger solver on uniform grid
for Gaussian contrast.
$N$ denotes the number of discretization points,
$e(x)=|u(x)-\tilde{u}(x)|$ 
is the error at the point $x$. $(0.5,0)$ is well inside the support of 
the contrast function and $(1,0.5)$ is a point where $q(x)$
has vanished to machine precision.
$Time$ is the time required by the fast solver. 
}
    \begin{tabular}{|c|r|r|r| }
    \hline
     $N$ & $e(0.5,0)$ & $e(1,0.5)$ & Time \\
    \hline\hline
    $256$  & $0.8925E+0$ & $0.4890E+0$ & $0.14$ \\
    $1024$ & $0.1809E+0$ & $0.1406E+0$ & $1.16$ \\
    $4096$ & $1.0938E-2$ & $9.6485E-3$ & $12.93$ \\    
    $16384$ & $3.1169E-4$ & $3.2633E-4$ & $104.1$  \\
    $65536$ & $1.4300E-5$ & $1.2537E-5$ & $844.8$  \\
    $262144$ & $8.8874E-7$ & $6.4485E-7$ & $6869.7$ \\
    \hline
    \end{tabular}\label{gauss_uniform}
    \end{center}
\end{table}

It is straightforward to check from Table
\ref{gauss_uniform} that the convergence rate
is approximately fourth order (for a target either within the support of 
$q(x)$ or in its exterior). The CPU requirements of the solver
can also be seen to scale approximately as $O(N^{3/2})$.

We now solve both the Gaussian and flat bump problems 
using our adaptive refinement strategy.
Tables \ref{gauss_adaptive} and \ref{erfc_adaptive} summarize the 
corresponding numerical results.

\begin{table}[h]
    \caption{Performance of adaptive solver for Gaussian contrast.
$\epsilon$ is the tolerance requested from adaptive discretization of the 
data. The other columns are the same as in Table \ref{gauss_uniform}.}
\begin{center}
\begin{tabular}{|c|r|r|r|r|r|r|r| }
\hline
$\epsilon$ & $N$ & $e(0.5,0)$ & $e(1,0.5)$ & Time \\
\hline\hline
$1E-4$ & $4096$  & $1.09E-2$ & $9.65E-3$ & $12.3$ \\
$1E-5$ & $4864$  & $7.89E-4$ & $1.20E-4$ & $16.2$ \\
$1E-6$ & $6400$  & $2.99E-4$ & $3.33E-4$ & $30.4$ \\    
$1E-7$ & $10240$ & $4.83E-5$ & $8.75E-6$ & $78.7$  \\
$1E-8$ & $16384$ & $1.13E-5$ & $5.35E-6$ & $190.1$ \\
$1E-9$ & $34816$ & $2.30E-6$ & $3.86E-7$ & $650.6$  \\
$1E-10$ & $70912$ & $7.14E-7$ & $3.16E-7$ & $2350.7$ \\                
$1E-11$ & $138688$ & $4.25E-8$ & $1.94E-8$ & $7363.4$ \\                
\hline
\end{tabular}\label{gauss_adaptive}
\end{center}
\end{table}

\begin{table}[h]
    \caption{Performance of adaptive solver for flat bump contrast.
$\epsilon$ is the tolerance requested from adaptive discretization of the 
data and $N$ is the total number of discretization points.
$e(x)$ is the error at the point $x$. $(0.5,2)$ is well inside the 
support of the contrast function and $(3,0.5)$ is a point where $q(x)$
has vanished to machine precision.}
\begin{center}
\begin{tabular}{|c|r|r|r|r|r|r|r| }
\hline
$\epsilon$ & $N$ & $e(0.5,2)$ & $e(3,0.5)$ & Time \\
\hline\hline
$1E-3$ & $10240$ & $1.88E-1$ & $9.29E-1$ & $61.8$ \\
$1E-4$ & $40192$   & $1.07E-2$ & $2.27E-2$ & $609.2$ \\    
$1E-5$ & $134656$ & $8.47E-4$ & $8.44E-4$ & $5561.7$ \\
$1E-7$ & $596224$ & $2.07E-5$ & $1.31E-5$ & $75861.5$ \\        
\hline
\end{tabular}\label{erfc_adaptive}
\end{center}
\end{table}

As expected, the adaptive method provides significant advantages:
for the same accuracy, the problem with a Gaussian contrast function 
is about ten times faster.

\begin{figure}[htbp]
\centering
\subfigure[Error at the point (0.5,2) as a function of $\kappa$.]{
\includegraphics[scale=.25]{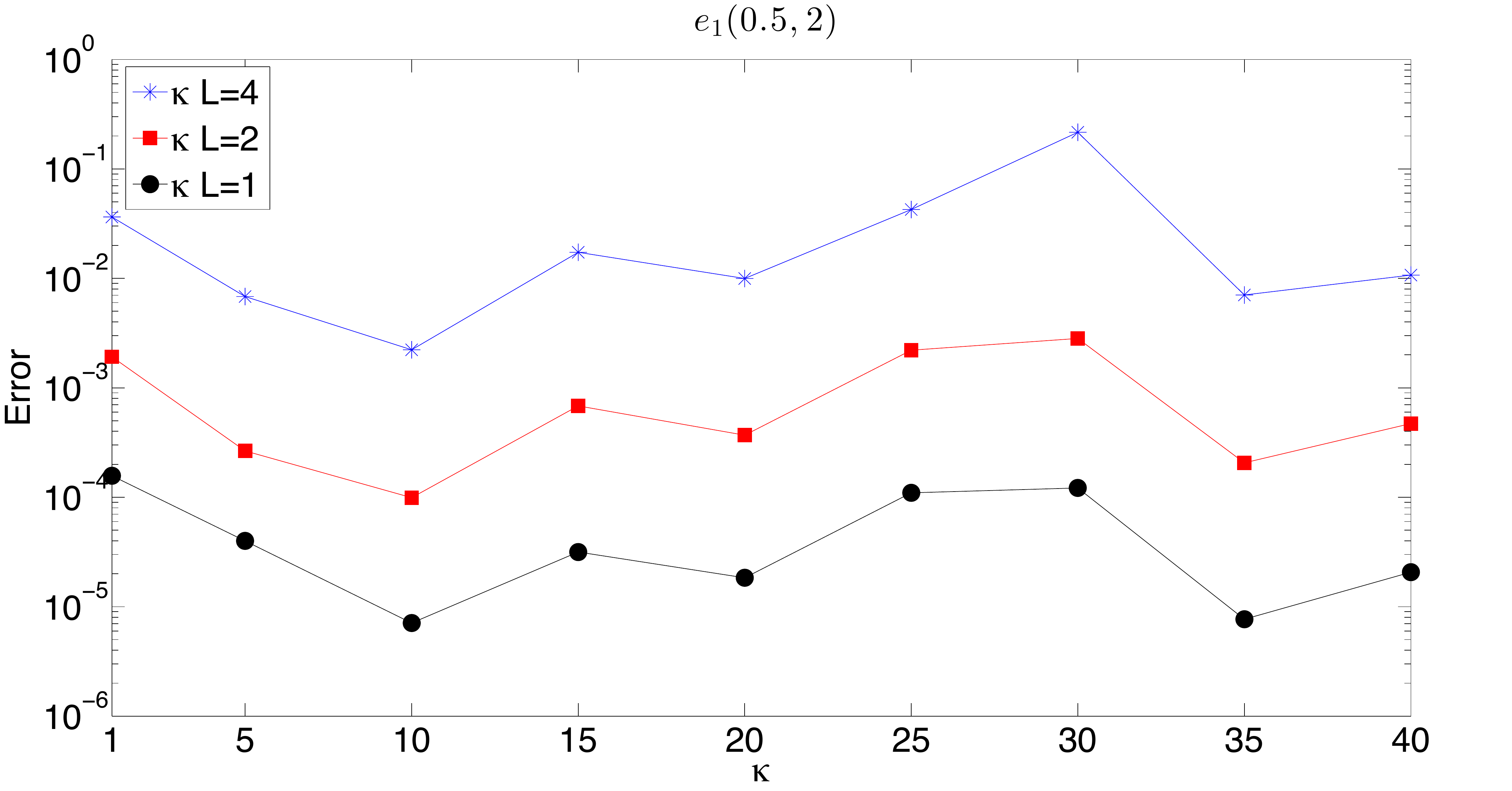}\label{figure_error_1_0}
}
\subfigure[Error at the point (3,0.5) as a function of $\kappa$.]{
\includegraphics[scale=.25]{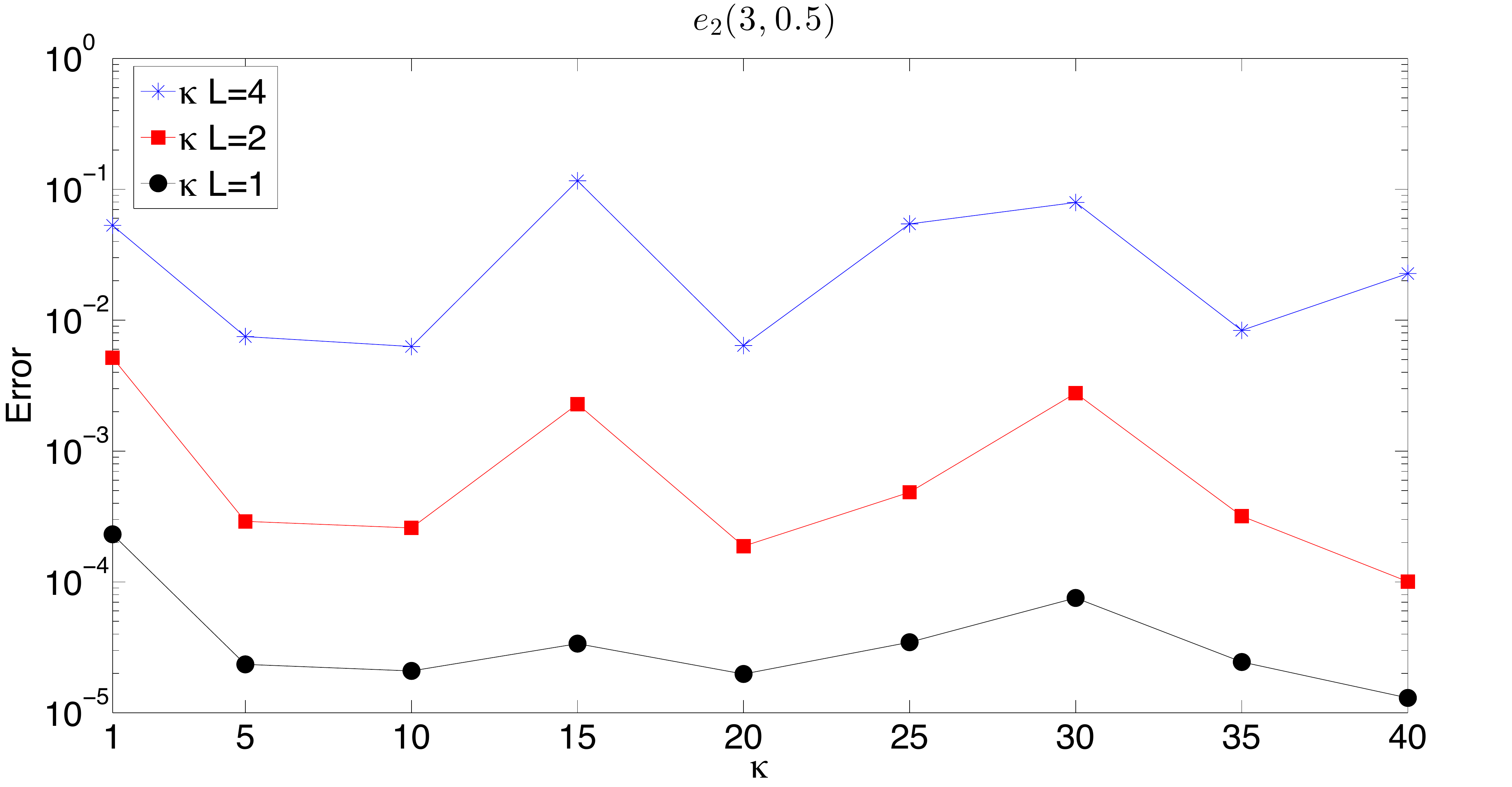}\label{figure_error_1_1}
}
\label{figure_error_1}
\caption{Error as a function of $\kappa$, which determines the size
of the domain in wavelengths. The values 
$\kappa\lb=4,2,1$ correspond to using approximately 6, 12 or 24 points per 
wavelength.}
\end{figure}

One of the advantages of using the Lippmann-Schwinger formulation
over direct discretization of the Helmholtz equation is that Green's function
based methods are not subject to the same kind of grid dispersion error.
To verify this, in 
Figs. \ref{figure_error_1_0} and \ref{figure_error_1_1}, we plot
the error for the flat bump contrast function as a function of 
wavenumber $\kappa$ at the points $(0.5,2)$ and $(3,0.5)$, respectively. 
Each curve represents the error for a fixed value of the parameter 
$\kappa\lb$ on leaf nodes. 
$\kappa\lb=4,2,1$ correspond to approximately 6, 12 and 24 points per 
wavelength. Note that the error does not grow with $\kappa$ as it would
with a fourth order discretization of the PDE.

\begin{example}Multiple scattering from well-separated Gaussian
bumps\end{example}

\begin{figure}[htbp]
\centering
\subfigure[Contrast function for twenty Gaussian bumps]{
\includegraphics[scale=.25]{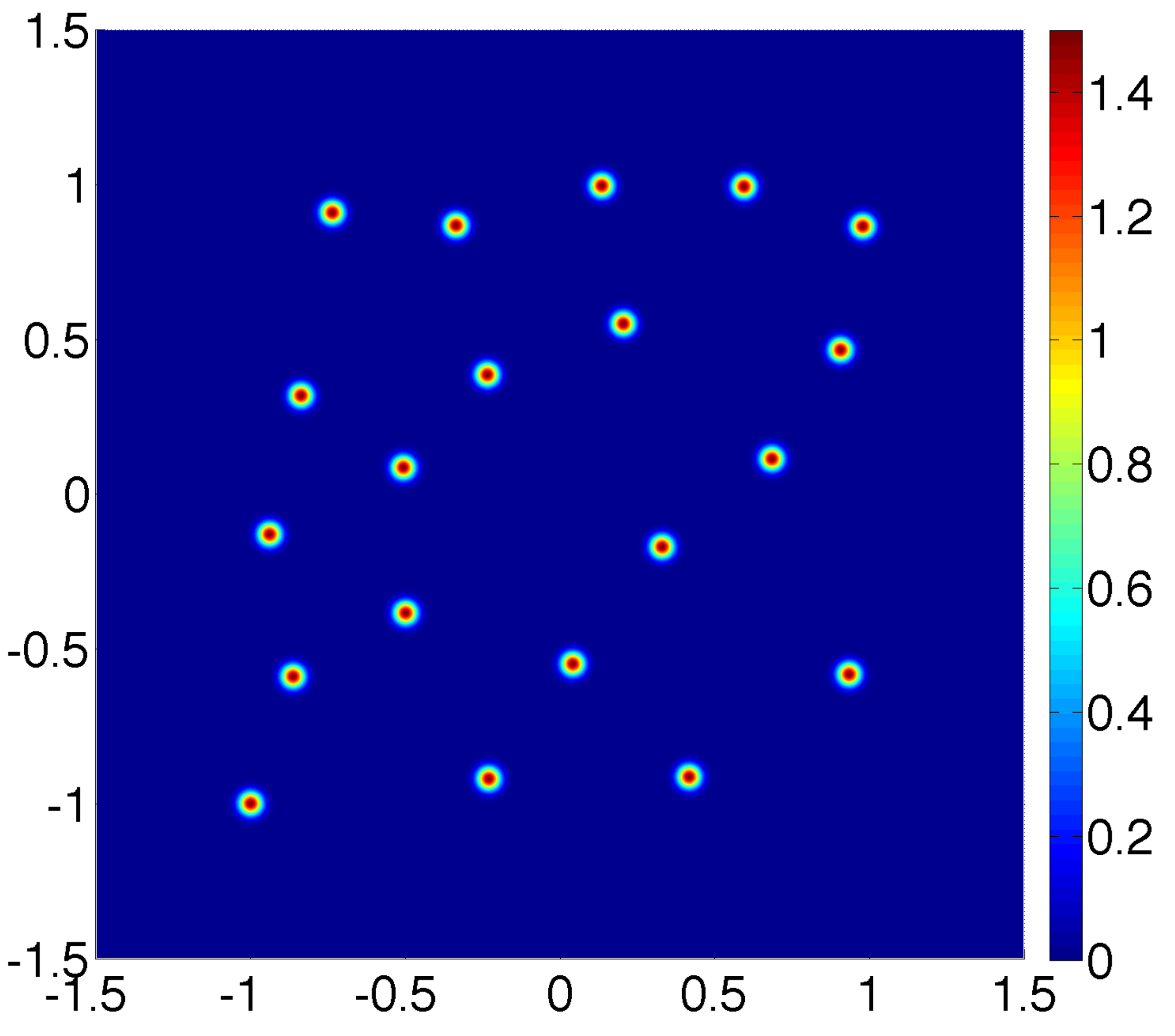}\label{figure_2_0}
}
\subfigure[Real part of total field]{
\includegraphics[scale=.25]{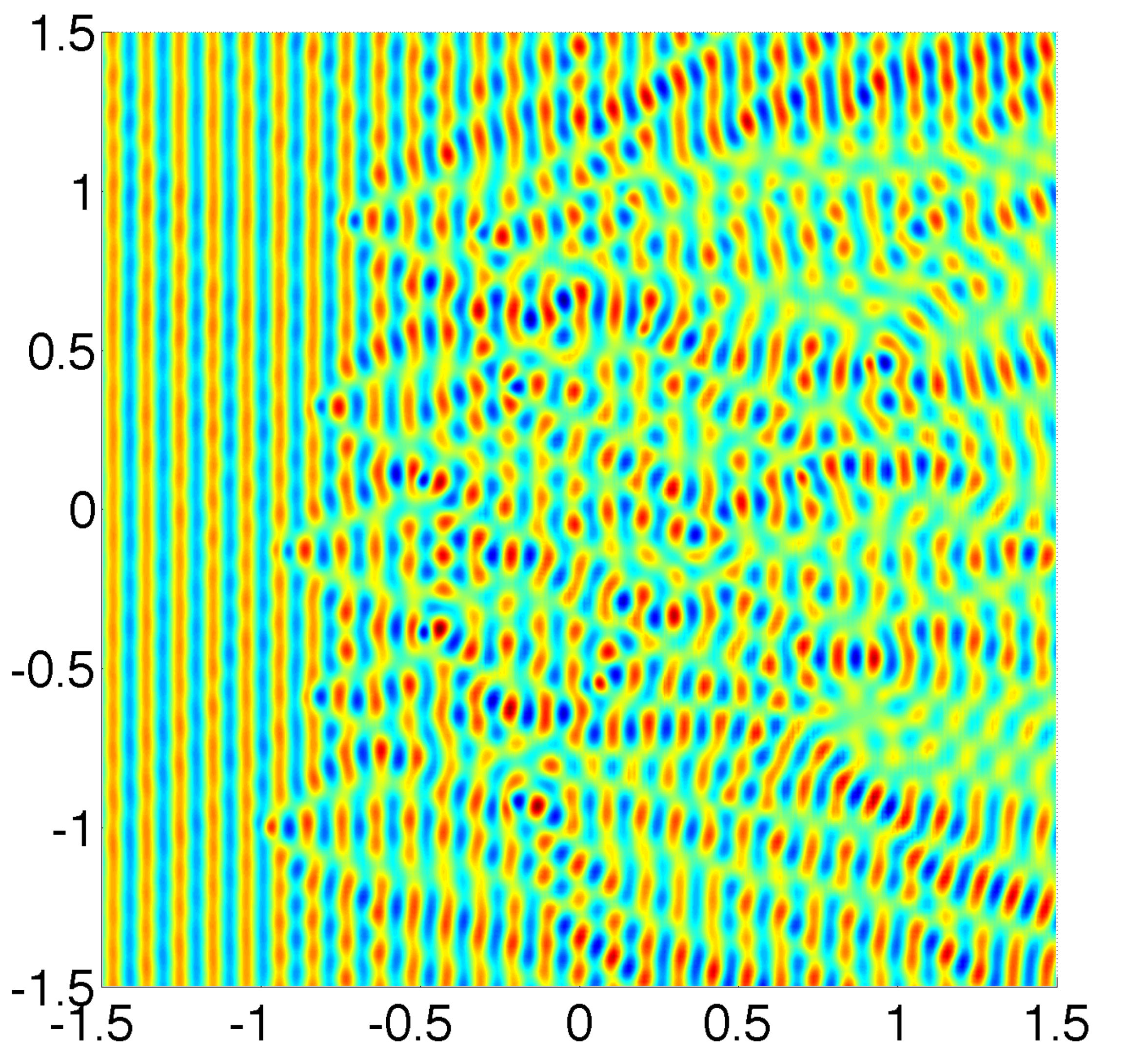}\label{figure_2_1}
}
\subfigure[Discretization of contrast function]{
\includegraphics[scale=1.75]{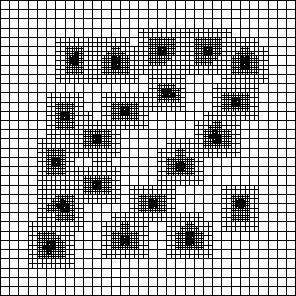}\label{figure_2_2}
}
\subfigure[Surface plot of contrast function]{
\includegraphics[scale=1]{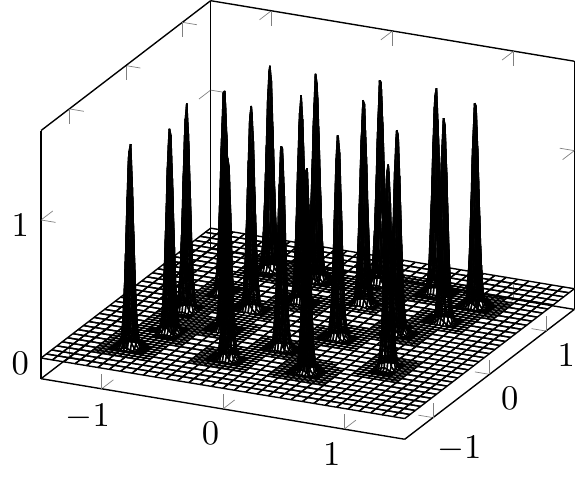}\label{figure_2_3}
}
\label{figure_example_2}
\caption{Multiple scattering from 20 Gaussian bumps}
\end{figure}

To illustrate the performance of the scheme on a more complicated 
example, we assume the contrast consists of the sum of $20$ Gaussian bumps 
of the form
 \begin{equation*}
q_j(x) = 1.5e^{-|x-x_j|^2/a},
 \end{equation*} 
where the centers $\{ x_j \}$ are randomly located in 
$[-1.5,1.5]^2$ and $a=0.0013$.
With wavenumber $\kappa=60$, the domain is approximately $30$ wavelengths in
each linear dimension. We discretize the domain using at least $9$ points 
per wavelength and resolve the contrast with tolerance $\epsilon$. 
In Fig. \ref{figure_2_0}, we plot the contrast function,
and in Fig. \ref{figure_2_1} we plot the real part of the total field 
calculated using $N=895264$ discretization points with about five digits 
of accuracy. Fig. \ref{figure_2_2} shows the adaptive grid (at a coarser
discretization) and 
Fig. \ref{figure_2_3} presents a surface plot of the contrast function. 
Table \ref{bumps_error} summarizes our numerical results. 

\begin{table}[htbp]
    \caption{Convergence behavior of scattering from multiple Gaussian
bumps.
$\epsilon$ is the tolerance requested from adaptive discretization of the 
data and $N$ is the number of discretization points.
$\Re(u(x))$ and 
$\Im(u(x))$ denote the real and imaginary parts of the solution at the
point $x$.}
\begin{center}
\begin{tabular}{|c|r|r|r|r|r| }
\hline
$\epsilon$ & $N$ & $\Re(u(0.5,2))$ & $\Im(u(0.5,2))$ & $\Re(u(3,0.5))$ & $\Im(u(3,0.5))$\\
\hline\hline
$1E-7$ & $155056$ & $-1.74310E-1$   & $-2.89411E-1$ & $-7.73067E-1$ & $-4.49033E-1$\\
$1E-8$ & $272128$   & $-1.74256E-1$ & $-2.89433E-1$ & $-7.73039E-1$ & $-4.49080E-1$\\
$1E-9$ & $441520$ & $-1.74268E-1$   & $-2.89425E-1$ & $-7.73045E-1$ & $-4.49084E-1$\\
$1E-10$ & $895264$ & $-1.74266E-1$   & $-2.89425E-1$ & $-7.73044E-1$ & $-4.49085E-1$\\
\hline
\end{tabular}\label{bumps_error}
\end{center}
\end{table}

\begin{example}An example from plasma physics\end{example}

Radio frequency waves play an important role in magnetic confinement 
fusion for both heating and diagnostic applications. 
Radio frequency wave propagation
in a fusion reactor can be described using a cold plasma model
\cite{Bonoli,Imbert-Gerard,Weitzner},
where the plasma is considered a fluid, possibly involving several species, 
subject to electric and magnetic forces. Here we consider only electrons, 
neglecting ions because of the mass ratio. 
Coupling Newton's law for the fluid motion with Maxwell's equations, driven
by a  current proportional to
the electrons velocity, leads to various propagation modes.
These are classified by their orientation with respect to the confining
magnetic field $\mathbf B_0$. The ordinary (O) mode corresponds to an
electric field parallel to $\mathbf B_0$, propagating in a plane orthogonal
to $\mathbf B_0$. This particular mode is used for reflectometry, 
a non-invasive method to probe the plasma density.

The O mode propagation in the cold plasma model reduces to the Helmholtz
equation \eqref{eqn_time_harmonic_Helm}, where the contrast $q(x)$ is
proportional to the additive inverse of the electron density, which is compactly
supported as the plasma is surrounded by vacuum.
Reflectometry relies on the fact that the incoming wave impinging on the
plasma can only propagate if the density is smaller than a threshold,
called the ``cut-off". In equation \eqref{eqn_time_harmonic_Helm},
the cut-off is implicitly defined by $q(x)=-1$. 
When the wave reaches the cut-off it
cannot propagate further and is reflected back from the plasma. In other 
words, the plasma forms an evanescent medium when the density is 
higher than the cut-off density, i.e. where $q(x)<-1$.

\begin{figure}[!htp]
\centering
\subfigure[Contrast function $q(x)$]{
\includegraphics[scale=.25]{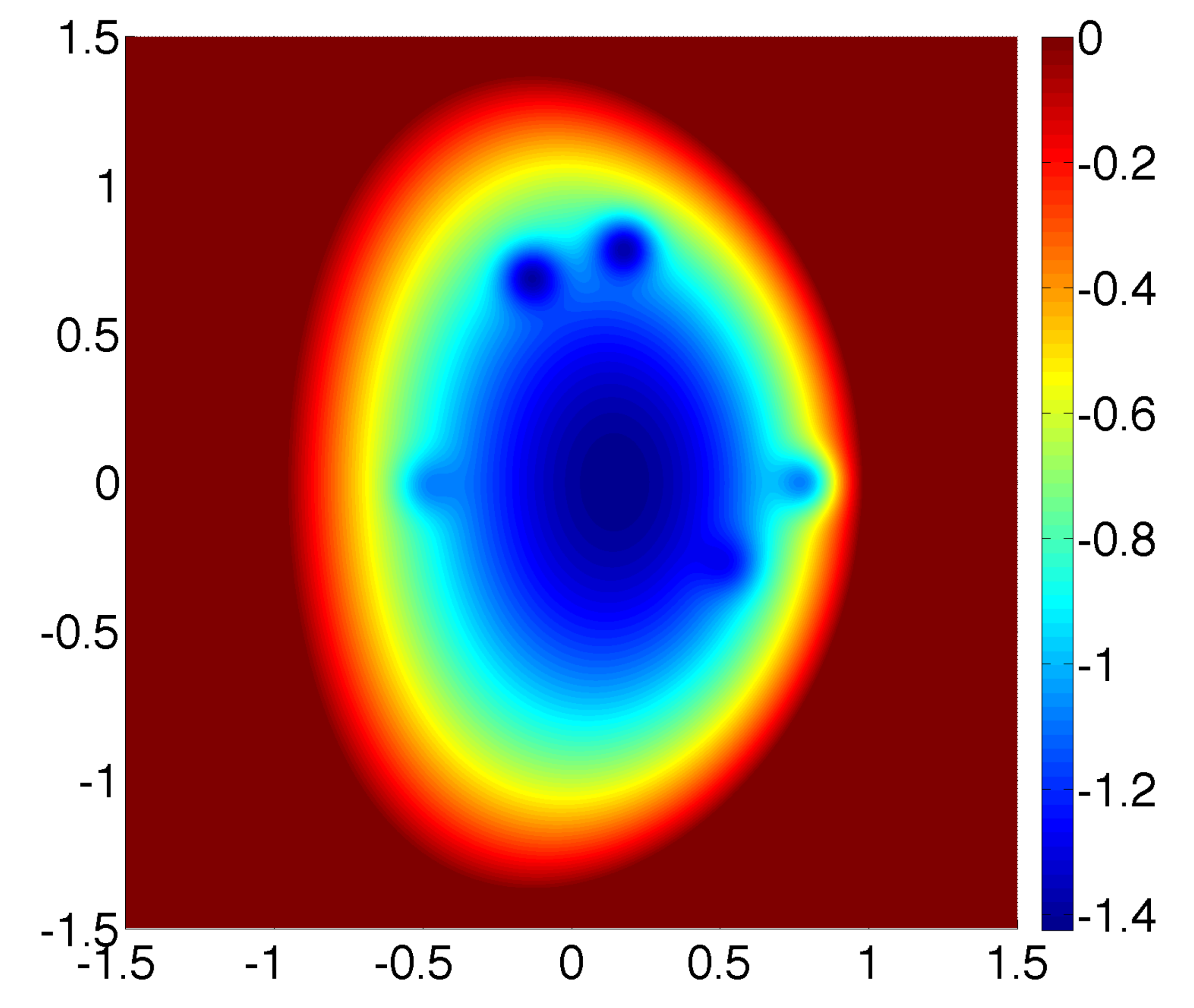}\label{figure_3_0}
}
\subfigure[Real part of total field ]{
\includegraphics[scale=.25]{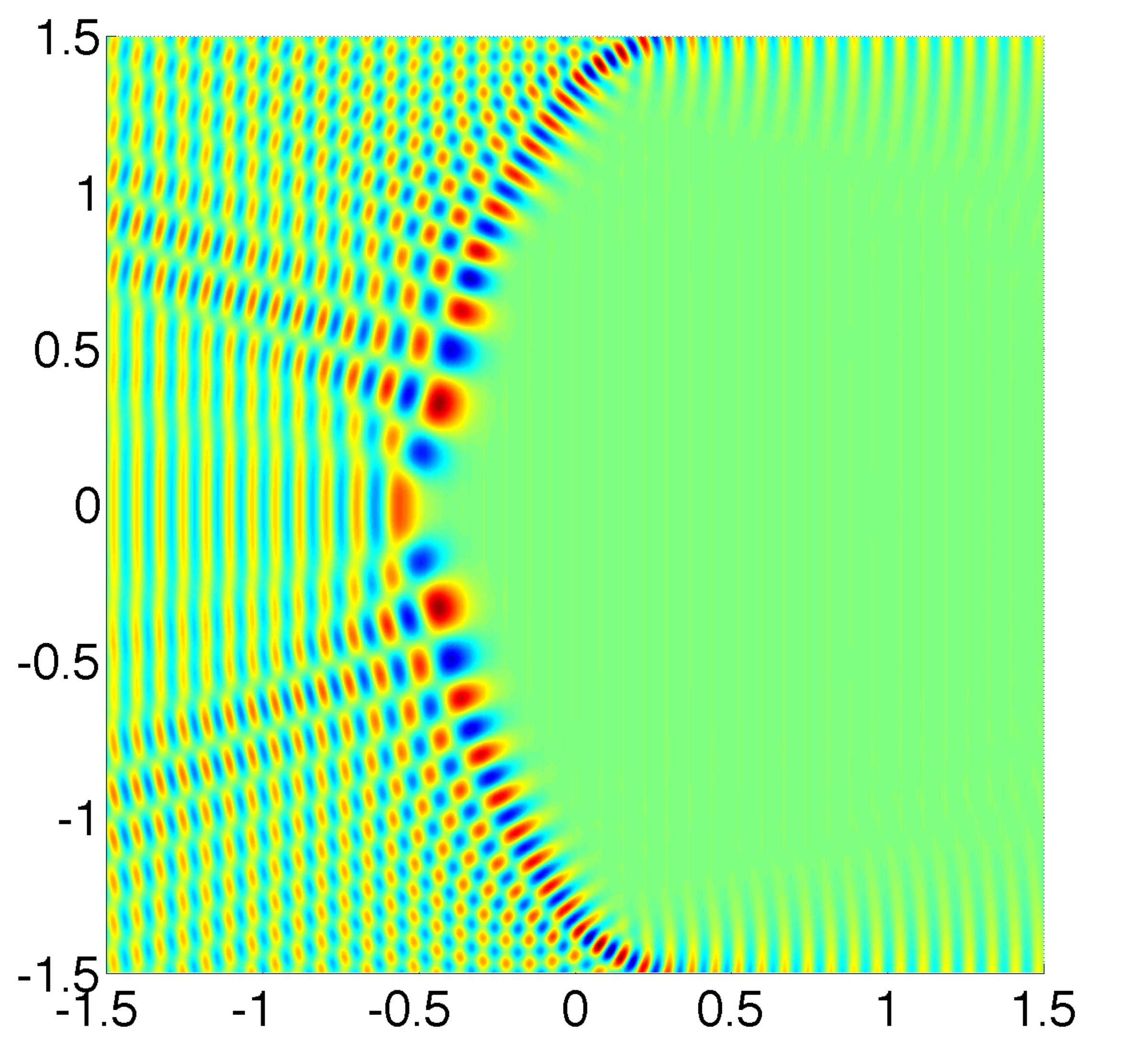}\label{figure_3_1}
}
\label{figure_example_3}
\caption{RF wave impinging on plasma.}
\end{figure}

In Fig. \ref{figure_3_0}, we plot the contrast for an idealized density
profile, defined by 
\begin{equation*}
q(x,y)=\twopartdef { 0, } {\psi(x,y) \leq 0.05,} {-1.5(\psi(x,y) - 0.05) - g(x,y)\cos(0.9y),} {\psi(x,y) > 0.05,}
%1.5(\psi(x,y) - 0.05) +\sum_{j=0}^5 a_j e^{-\frac{((x-x_j)^2+(y-y_j))}{0.01}}\cos(0.9y)
\end{equation*}
where 
\begin{align*}
\psi(x,y) & = 1-(x-0.15(1-x^2))^2 - C(1+0.3x)^2 y^2, \\
g(x,y)     & = \sum_{j=1}^5 a_j \exp\left(-\frac{(x-x_j)^2+(y-y_j)^2}{0.01}\right),
 \end{align*}
with $C=0.4987$, and the values of $a_j$, $x_j$, and $y_j$, for $j=1,\ldots,5$, are given in Table \ref{gauss_tab}.
We set the wavenumber $\kappa=85$. To solve this problem, we used $N=1060864$ 
points corresponding to approximately $25$ points per wavelength. 
The total field calculated by our solver is
shown in Figure \ref{figure_3_1}. As expected,  there is no
propagation through the evanescent zone and the incident wave is reflected
back at the cut-off.
Convergence analysis using Richardson extrapolation indicates that
we obtained more than 6 digits of
accuracy. 

\begin{table}[htbp]
    \caption{Parameters $a_j$, $x_j$, and $y_j$ for the contrast function used in the plasma physics example.}
\begin{center}
\begin{tabular}{|c|r|r|r|r|r| }
\hline
$j$ & $a_j$ & $x_j$ & $y_j$ \\
\hline\hline
$1$ & $0.45$   & $0.80$  & $0.00$\\
$2$ & $0.195$ & $0.54$  & $-0.28$\\
$3$ & $0.51$   & $-0.14$ & $0.70$  \\
$4$ & $0.195$ & $-0.50$ & $-0.01$\\
$5$ & $0.63$   & $0.18$  & $0.80$\\
\hline
\end{tabular}\label{gauss_tab}
\end{center}
\end{table}

\section{Conclusion} \label{sec:concl}

We have presented a high-order method for the 
Lippmann-Schwinger equation that is of particular value when coupled 
with modern fast, direct solvers. More precisely, we have developed
algorithms for rapidly accessing arbitrary elements of the system
matrix consistent with an adaptive Nystr\"{o}m discretization. 
Our approach extends naturally to
three dimensions and to other governing Green's functions, although
some of the accelerations have been developed specifically for the 
Helmholtz equation and will require modification in other settings.

Our numerical examples were computed using a fourth order accurate approach, 
but the extension to arbitrary order is straightforward.
We are currently working on 
the coupling of our tools with fast solvers that are quad-tree
or oct-tree based, rather than k-d tree-based as in our HODLR solver.
These were discussed very briefly in section
\ref{section_fast_direct}.
This will substantially reduce the constant implicit in the 
$O(N^{3/2})$ notation. 
We are also working on fast discretization of the 
Lippmann-Schwinger equation in three dimensions for both acoustic and
electromagnetic applications.

\section*{Acknowledgments}
This work was supported in part by the 
Applied Mathematical Sciences Program of the U.S. Department of Energy
under Contract 
DEFGO288ER25053 and
by the Office of the Assistant Secretary of Defense for Research and 
Engineering and AFOSR under NSSEFF Program Award FA9550-10-1-0180.
The authors would like to thank 
Michael O'Neil and Alex Barnett for several useful conversations.
\bibliography{./Bibnew.bib}

\begin{thebibliography}{10}

\bibitem{abramowitz}
{\sc M.~Abramowitz and I.~Stegun}, {\em Handbook of Mathematical Functions with
  Formulas, Graphs, and Mathematical Tables}, U.S. Government Printing Office,
  1964.

\bibitem{aguilar2002}
{\sc J.~Aguilar and Y.~Chen}, {\em High-order corrected trapezoidal quadrature
  ruless for functions with a logarithmic singularity in 2-d}, Computers \&
  Mathematics with Applications, 44 (2002), pp.~1031--1039.

\bibitem{Amar1983}
{\sc M.~B. Amar and F.~C. Farnoux}, {\em Numerical solution of the
  lippmann-schwinger equations in photoemission: application to xenon}, Journal
  of Physics B: Atomic and Molecular Physics, 16 (1983), p.~2339.

\bibitem{siva}
{\sc S.~Ambikasaran and E.~Darve}, {\em An ${O} ({N} log {N})$ fast direct
  solver for partial hierarchically semi-separable matrices --- with
  application to radial basis function interpolation.}, J. Sci. Comput.,
  (2013), pp.~477--501.

\bibitem{amirhossein2014fast}
{\sc A.~Aminfar, S.~Ambikasaran, and E.~Darve}, {\em A fast block low-rank
  dense solver with applications to finite-element matrices}, arXiv preprint
  arXiv:1403.5337 [cs-NA],  (2014).

\bibitem{bebendorf2005hierarchical}
{\sc M.~Bebendorf}, {\em Hierarchical {LU} decomposition-based preconditioners
  for {BEM}}, Computing, 74 (2005), pp.~225--247.

\bibitem{Bonoli}
{\sc P.~T. Bonoli}, {\em Electromagneti mode conversion: understanding waves
  that suddenly change thier nature}, Journal of Physics: Conference Series, 16
  (2005), pp.~35--39.

\bibitem{borm2003hierarchical}
{\sc S.~B{\"o}rm, L.~Grasedyck, and W.~Hackbusch}, {\em Hierarchical matrices},
  Lecture notes, 21 (2003).

\bibitem{borm2003introduction}
\leavevmode\vrule height 2pt depth -1.6pt width 23pt, {\em Introduction to
  hierarchical matrices with applications}, Engineering Analysis with Boundary
  Elements, 27 (2003), pp.~405--422.

\bibitem{Bruno2003}
{\sc O.~P. Bruno}, {\em Wave scattering by inhomogeneous media: efficient
  algorithms and applications}, Physica B: Condensed Matter, 338 (2003), pp.~67
  -- 73.
\newblock Proceedings of the Sixth International Conference on Electrical
  Transport and Optical Properties of Inhomogeneous Media.

\bibitem{Bruno2004}
{\sc O.~P. Bruno and E.~Hyde}, {\em An efficient, preconditioned, high-order
  solver for scattering by two-dimensional inhomogeneous media}, Journal of
  Computational Physics, 200 (2004), pp.~670 -- 694.

\bibitem{Bruno2005}
{\sc O.~P. Bruno and E.~M. Hyde}, {\em Higher-order fourier approximation in
  scattering by two-dimensional, inhomogeneous media.}, SIAM J. Numerical
  Analysis, 42 (2005), pp.~2298--2319.

\bibitem{chandrasekaran2006fast1}
{\sc S.~Chandrasekaran, P.~Dewilde, M.~Gu, W.~Lyons, and T.~Pals}, {\em A fast
  solver for {HSS} representations via sparse matrices}, SIAM Journal on Matrix
  Analysis and Applications, 29 (2006), pp.~67--81.

\bibitem{chandrasekaran2006fast}
{\sc S.~Chandrasekaran, M.~Gu, and T.~Pals}, {\em A fast {ULV} decomposition
  solver for hierarchically semiseparable representations}, SIAM Journal on
  Matrix Analysis and Applications, 28 (2006), pp.~603--622.

\bibitem{chen2002fast}
{\sc Y.~Chen}, {\em A fast, direct algorithm for the lippmann--schwinger
  integral equation in two dimensions}, Advances in Computational Mathematics,
  16 (2002), pp.~175--190.

\bibitem{huang2006}
{\sc H.~Cheng, J.~Huang, and T.~J. Leiterman}, {\em An adaptive fast solver for
  the modified helmholtz equation in two dimensions}, Journal of Computational
  Physics, 211 (2006), pp.~616--637.

\bibitem{Colton}
{\sc D.~Colton and R.~Kress}, {\em Inverse Acoustic and Electromagnetic
  Scattering Theory}, Springer, 2$^\text{nd}$~ed., 1998.

\bibitem{corona2014n}
{\sc E.~Corona, P.-G. Martinsson, and D.~Zorin}, {\em An o (n) direct solver
  for integral equations on the plane}, Applied and Computational Harmonic
  Analysis,  (2014).

\bibitem{Costabel}
{\sc M.~Costabel}, {\em On the spectrum of volume integral operators in
  acoustic scattering}, in Integral Methods in Science and Engineering,
  Birkha\"{u}ser, Chapter 11, 2015.

\bibitem{crutchfield2006remarks}
{\sc W.~Crutchfield, Z.~Gimbutas, L.~Greengard, J.~Huang, V.~Rokhlin,
  N.~Yarvin, and J.~Zhao}, {\em Remarks on the implementation of wideband fmm
  for the helmholtz equation in two dimensions}, Contemporary Mathematics, 408
  (2006), pp.~99--110.

\bibitem{duan2009}
{\sc R.~Duan and V.~Rokhlin}, {\em High order quadratures for the solution of
  scattering problems in two dimensions}, Journal of Computational Physics, 228
  (2009), pp.~2152--2174.

\bibitem{ey11}
{\sc B.~Engquist and L.~Ying}, {\em Sweeping preconditioner for the helmholtz
  equation: hierarchical matrix representation}, Communications on Pure and
  Applied Mathematics, 64 (2011), pp.~697--735.

\bibitem{Ethridge2001}
{\sc F.~Ethridge and L.~Greengard}, {\em {A new fast-multipole accelerated
  Poisson solver in two dimensions}}, SIAM Journal on Scientific Computing, 23
  (2001), pp.~741--760.

\bibitem{Barnett2015}
{\sc A.~Gillman, A.~H. Barnett, and P.-G. Martinsson}, {\em A spectrally
  accurate direct solution technique for frequency-domain scattering problems
  with variable media}, BIT Numerical Mathematics, 55 (2015), pp.~141--170.

\bibitem{goreinov1997theory}
{\sc S.~Goreinov, E.~Tyrtyshnikov, and N.~Zamarashkin}, {\em A theory of
  pseudoskeleton approximations}, Linear Algebra and its Applications, 261
  (1997), pp.~1--21.

\bibitem{greengard2009fast}
{\sc L.~Greengard, D.~Gueyffier, P.-G. Martinsson, and V.~Rokhlin}, {\em Fast
  direct solvers for integral equations in complex three-dimensional domains},
  Acta {N}umerica, 18 (2009), pp.~243--275.

\bibitem{hackbusch2001introduction}
{\sc W.~Hackbusch, L.~Grasedyck, and S.~B{\"o}rm}, {\em An introduction to
  hierarchical matrices}, Max-Planck-Inst. f{\"u}r Mathematik in den
  Naturwiss., 2001.

\bibitem{hoying1}
{\sc K.~Ho and L.~Ying}, {\em Hierarchical interpolative factorization for
  elliptic operators: differential equations}, Communications in Pure and
  Applied Mathematics,  (to appear).

\bibitem{hoying2}
\leavevmode\vrule height 2pt depth -1.6pt width 23pt, {\em Hierarchical
  interpolative factorization for elliptic operators: integral equations},
  Communications in Pure and Applied Mathematics,  (to appear).

\bibitem{ho2012fast}
{\sc K.~L. Ho and L.~Greengard}, {\em A fast direct solver for structured
  linear systems by recursive skeletonization}, SIAM Journal on Scientific
  Computing, 34 (2012), pp.~2507--2532.

\bibitem{Hohage2006}
{\sc T.~Hohage}, {\em {Fast numerical solution of the electromagnetic medium
  scattering problem and applications to the inverse problem}}, Journal of
  Computational Physics, 214 (2006), pp.~224--238.

\bibitem{Imbert-Gerard}
{\sc L.-M. Imbert-Gerard}, {\em Well-posedness of 2d mode conversion model and
  generalized plaen wave simulations}, Journal of Physics: Conference Series,
  (preprint).

\bibitem{kong2011adaptive}
{\sc W.~Y. Kong, J.~Bremer, and V.~Rokhlin}, {\em An adaptive fast direct
  solver for boundary integral equations in two dimensions}, {A}pplied and
  {C}omputational {H}armonic {A}nalysis, 31 (2011), pp.~346--369.

\bibitem{Lanzara2004}
{\sc F.~Lanzara, V.~G. Maz´ya, and G.~Schmidt}, {\em Numerical solution of the
  lippmann-schwinger equation by approximate approximations}, Journal of
  Fourier Analysis and Applications, 10 (2004), pp.~645--660.

\bibitem{liberty2007randomized}
{\sc E.~Liberty, F.~Woolfe, P.-G. Martinsson, V.~Rokhlin, and M.~Tygert}, {\em
  Randomized algorithms for the low-rank approximation of matrices},
  Proceedings of the National Academy of Sciences, 104 (2007), p.~20167.

\bibitem{Biros}
{\sc D.~Malhotra, A.~Gholami, and G.~Biros}, {\em A volume integral equation
  stokes solver for problems with variable coefficients}, in International
  Conference for High Performance Computing, Networking, Storage and Analysis,
  SC14, IEEE, Nov. 2014, pp.~92--102.

\bibitem{martinsson2009fast}
{\sc P.-G. Martinsson}, {\em A fast direct solver for a class of elliptic
  partial differential equations}, Journal of Scientific Computing, 38 (2009),
  pp.~316--330.

\bibitem{martinsson2013}
\leavevmode\vrule height 2pt depth -1.6pt width 23pt, {\em A direct solver for
  variable coefficient elliptic pdes discretized via a composite spectral
  collocation method}, Journal of Computational Physics, 242 (2013),
  pp.~460--479.

\bibitem{martinsson2005fast}
{\sc P.-G. Martinsson and V.~Rokhlin}, {\em A fast direct solver for boundary
  integral equations in two dimensions}, Journal of Computational Physics, 205
  (2005), pp.~1--23.

\bibitem{handbook}
{\sc F.~W.~J. Olver, D.~W. Lozier, R.~F. Boisvert, and e.~C.~W.~Clark}, {\em
  NIST Handbook of Mathematical Functions}, Cambridge University Press, New
  York, 2010.

\bibitem{rjasanow2002adaptive}
{\sc S.~Rjasanow}, {\em Adaptive cross approximation of dense matrices}, in
  Int. Association Boundary Element Methods Conf., IABEM, 2002, pp.~28--30.

\bibitem{Sandfort2010}
{\sc K.~Sandfort}, {\em The factorization method for inverse scattering from
  periodic inhomogeneous media}, KIT Scientific Publishing, Karlsruhe, 2010.

\bibitem{vainikko2000fast}
{\sc G.~Vainikko}, {\em Fast solvers of the lippmann-schwinger equation}, in
  Direct and inverse problems of mathematical physics, Springer, 2000,
  pp.~423--440.

\bibitem{Weitzner}
{\sc H.~Weitzner}, {\em Lower hybrid waves in the cold plasma model},
  Communications on Pure and Applied Mathematics, 38 (1985), pp.~919--932.

\bibitem{zd14}
{\sc L.~Zepeda-Nunez and L.~Demanet}, {\em The method of polarized traces for
  the 2d helmholtz equation}, submitted,  (2014), pp.~1--54.

\bibitem{zhao2005adaptive}
{\sc K.~Zhao, M.~N. Vouvakis, and J.-F. Lee}, {\em The adaptive cross
  approximation algorithm for accelerated method of moments computations of emc
  problems}, Electromagnetic Compatibility, IEEE Transactions on, 47 (2005),
  pp.~763--773.

\end{thebibliography}

\end{document}